\documentclass[12pt]{article}
\usepackage{amssymb}
\usepackage{amsmath}
\usepackage{epsfig}
\usepackage{theorem}

\newtheorem{thm}{Theorem}[section]
\newtheorem{prop}[thm]{Proposition}

\newtheorem{prob}[thm]{Problem}

{\theorembodyfont{\rm}
\newtheorem{cor}[thm]{Corollary}

\newtheorem{rem}[thm]{Remark}}

\newcommand{\ra}{\rightarrow}
\newcommand{\supp}{\mathrm{supp}}

\newcommand{\PGL}{\mathrm{PGL}}
\newcommand{\Isom}{\mathrm{Isom}}

\newcommand{\Hom}{\mathrm{Hom}}
\newcommand{\Ext}{\mathrm{Ext}}
\newcommand{\Pic}{\mathrm{Pic}}
\newcommand{\Spec}{\mathrm{Spec\ }}
\newcommand{\Proj}{\mathrm{Proj}}
\newcommand{\Sym}{\mathrm{Sym}}

\newcommand{\Gr}{\mathrm{Gr}}
\newcommand{\SL}{\mathrm{SL}}

\newcommand{\bC}{{\mathbb{C}}}
\newcommand{\bG}{{\mathbb{G}}}
\newcommand{\bL}{{\mathbb{L}}}
\newcommand{\bP}{{\mathbb{P}}}
\newcommand{\bQ}{{\mathbb{Q}}}
\newcommand{\bR}{{\mathbb{R}}}
\newcommand{\bZ}{{\mathbb{Z}}}
\newcommand{\cA}{{\mathcal{A}}}
\newcommand{\cB}{{\mathcal{B}}}
\newcommand{\cC}{{\mathcal{C}}}
\newcommand{\cD}{{\mathcal{D}}}

\newcommand{\cL}{{\mathcal{L}}}
\newcommand{\cM}{{\mathcal{M}}}

\newcommand{\cO}{{\mathcal{O}}}
\newcommand{\cQ}{{\mathcal{Q}}}
\newcommand{\cT}{{\mathcal{T}}}
\newcommand{\cW}{{\mathcal{W}}}
\newcommand{\hC}{{\hat C}}
\newcommand{\hF}{{\hat F}}
\newcommand{\hT}{{\hat T}}
\newcommand{\hZ}{{\hat Z}}

\newcommand{\ocM}{\overline {\mathcal{M}}}
\newcommand{\oL}{\overline L}
\newcommand{\oM}{\overline M}

\author{Brendan Hassett
\thanks{
The author was supported by the Institute
of Mathematical Sciences of the Chinese University of Hong
Kong, an NSF Postdoctoral Fellowship, and NSF Continuing 
Grant 0070537. 
}}

\title{Moduli spaces of weighted pointed stable curves}

\date{March 2, 2001\\
Revised March 14, 2002}

\begin{document}
\maketitle

\section{Introduction}

It has long been understood that a moduli space may admit a plethora
of different compactifications, each corresponding to a choice 
of combinatorial data.  Two outstanding examples are the
toroidal compactifications of quotients of bounded symmetric
domains \cite{AMRT} and the theory of variation of geometric invariant 
theory (GIT) quotients \cite{BP} \cite{DH} \cite{Th}.  
However, in both of these situations
a modular interpretation of the points added at the boundary can be
elusive.  By a modular interpretation, we mean the description of
a moduli functor whose points are represented by the compactification.
Such moduli functors should naturally incorporate the combinatorial
data associated with the compactification.

The purpose of this paper is to explore in depth one case where
functorial interpretations are readily available:  configurations
of nonsingular points on a curve.  Our standpoint is to consider
pointed curves as `log varieties', pairs $(X,D)$ where $X$ is 
a variety and $D=\sum_i a_iD_i$ is an effective ${\bQ}$-divisor on $X$.  
The minimal model program suggests a construction for the moduli space of
such pairs provided they are {\em stable}, i.e., $(X,D)$
should have relatively mild singularities and the divisor $K_X+D$
should be ample.  Of course, when $\dim(X)=1$ and $D$ is reduced
the resulting moduli space is the Mumford-Knudsen moduli space
of pointed stable curves \cite{KnMu}\cite{Kn}\cite{Kn2}.  
In \S \ref{sect:construct} we give a construction 
for arbitrary $D$.  When $\dim(X)=2$,
a proof for the existence of such moduli spaces was given by Koll\'ar, 
Shepherd-Barron, and Alexeev \cite{Ko1},\cite{Al1} when $D=0$,
and was sketched by Alexeev\cite{Al2} for $D\ne 0$.  The case
of higher dimensions is still open, but would follow
from standard conjectures of the minimal model program \cite{Karu}.   

We embark on a systematic study of the dependence of these 
moduli spaces on the coefficients of the divisor $D$.  We find natural
transformations among the various moduli functors which induce
birational reduction morphisms among the associated 
compactifications (see \S \ref{sect:reduction}).  These
morphisms can often be made very explicit.  We recover the
alternate compactifications studied by Kapranov \cite{Kap1} \cite{Kap2}, 
Keel\cite{Ke}, and Losev-Manin \cite{LM} as special cases of our 
theory (see \S \ref{sect:examples}).  The blow-up constructions
they describe are closely intertwined with our functorial reduction maps.  
The resulting contractions may sometimes
be understood as log minimal models of the moduli space itself,
where the log divisor is supported in the boundary (see \S \ref{sect:LMMP}).

The moduli spaces we consider do not obviously
admit a uniform construction as the quotients arising from varying 
the linearization of an invariant theory problem.  However,
ideas of Kapranov (see \cite{Kap1} 0.4.10) suggest
indirect GIT approaches to our spaces.  Furthermore, we indicate
how certain GIT quotients may be interpretted as `small parameter limits' of our
moduli spaces, and the flips between these GIT quotients factor naturally
through our spaces (see \S \ref{sect:GIT}).     

One motivation for this work is the desire for a better understanding
of compactifications of moduli spaces of log surfaces.  These
have been studied in special cases \cite{Has} and it was found
that the moduli space depends on the coefficients of the boundary
in a complicated way.  For example, in the case of quintic
plane curves (i.e., $X={\bP}^2$ and $D=aC$ with $C$ a plane quintic) even
the irreducible component structure and dimension of the moduli space  
depends on $a$.  For special values of $a$ 
the moduli space sprouts superfluous irreducible components
attached at infinity.  This pathology is avoided when the
coefficient is chosen generically.  Furthermore, recent exciting work 
of Hacking \cite{Hac} shows that for small values of the coefficient $a$
the moduli space is often nonsingular and its boundary admits an explicit description.
Roughly, Hacking considers the moduli space parametrizing pairs $({\bP}^2,a C)$
where $C$ is a plane curve of degree $d$ as $a \ra 3/d$.  In a future paper, 
we shall consider birational transformations of moduli
spaces of log surfaces induced by varying $a$.

\     

\noindent {\bf Acknowledgments:}  
I would like to thank M. Kapranov for patiently discussing
his work, Paul Hacking for sharing an early version of
his thesis, 
Alessio Corti for kindly providing an informative resum\'e
of Hacking's work, and Igor Dolgachev for many helpful comments,
especially on \S \ref{sect:GIT}.  
I would also like to thank 
Linda Chen, S\'andor Kov\'acs, Yuri Tschinkel, Eckart Viehweg,
and Kang Zuo for helpful conversations.  This work greatly
benefitted from correspondence with S. Keel.  His ideas concerning 
the log minimal model program on moduli spaces
informs this work;  he also pointed out errors in 
an early version of \S \ref{sect:LMMP}.

\section{The moduli problem}

Fix nonnegative integers $g$ and $n$ and
let $B$ be a noetherian scheme. 
A family of nodal curves
of genus $g$ with $n$ marked points over $B$
consists of 
\begin{enumerate}
\item{a flat proper morphism $\pi:C \ra B$
whose geometric fibers are nodal connected 
curves of arithmetic genus $g$; and}
\item{sections $s_1,s_2,\ldots,s_n$ of $\pi$.}
\end{enumerate}
A morphism of two such families
$$\phi:(C,s_1,\ldots,s_n)\longrightarrow
(C',s'_1,\ldots,s'_n)$$ 
consists of a $B$-morphism $\phi:C \ra C'$ 
such that $\phi(s_j)\subset s'_j$ for
$j=1,\ldots,n$.  The set of isomorphisms of two
such families is denoted 
$$\Isom((C,s_1,\ldots,s_n),(C',s'_1,\ldots,s'_n)),$$
or simply $\Isom(C,C')$ when there is no risk of confusion.  

A collection of input data $(g,\cA):=(g,a_1,\ldots,a_n)$ consists
of an integer $g\ge 0$ and the {\em weight data},
an element $(a_1,\ldots,a_n)\in \bQ^n$
such that $0<a_j\le 1$ for $j=1,\ldots,n$ and 
$$2g-2+a_1+a_2+\ldots +a_n>0.$$
A family of nodal curves with marked points
$(C,s_1,\ldots,s_n)\stackrel{\pi}{\ra}B$
is stable of type $(g,\cA)$ if 
\begin{enumerate}
\item{the sections $s_1,\ldots,s_n$ lie in the smooth locus
of $\pi$, and for any subset $\{ s_{i_1},\ldots, s_{i_r}\}$ with nonempty
intersection we have $a_{i_1}+\ldots +a_{i_r}\le 1$;}
\item{$K_{\pi}+a_1s_1+a_2s_2+\ldots +a_ns_n$ is $\pi$-relatively
ample.}
\end{enumerate}
This coincides with the traditional notion of pointed
stable curves when $a_1=a_2=\ldots=a_n=1$.

\begin{thm}\label{thm:existence}
Let $(g,\cA)$ be a collection of input data.
There exists a connected Deligne-Mumford stack 
$\ocM_{g,\cA}$, smooth and proper over $\bZ$, 
representing the moduli problem of pointed
stable curves of type $(g,\cA)$.  The corresponding coarse
moduli scheme $\oM_{g,\cA}$ is projective over $\bZ$.  
\end{thm}
The universal curve is denoted 
$\cC_{g,\cA}\ra \ocM_{g,\cA}$.  
Theorem \ref{thm:existence} is proved in \S \ref{sect:construct}.  

\subsection{Variations on the moduli problem}
\subsubsection{Zero weights}\label{subsect:zerowts}
One natural variant on our moduli problem is to allow some
of the sections to have weight zero.
We consider $(g,{\tilde \cA}):=(g,a_1,\ldots,a_n)$ where
$(a_1,\ldots,a_n)\in \bQ^n$
with $0\le a_j\le 1$ and 
$$2g-2+a_1+a_2+\ldots +a_n>0.$$
A family of nodal curves with marked points
$(C,s_1,\ldots,s_n)\stackrel{\pi}{\ra}B$
is stable of type $(g,{\tilde \cA})$ if 
\begin{enumerate}
\item{the sections $s_i$ with positive weights lie in the smooth locus
of $\pi$, and for any subset $\{ s_{i_1},\ldots, s_{i_r}\}$ with nonempty
intersection we have $a_{i_1}+\ldots +a_{i_r}\le 1$;}
\item{$K_{\pi}+a_1s_1+a_2s_2+\ldots +a_ns_n$ is $\pi$-relatively
ample.}
\end{enumerate}
There is no difficulty making sense of the divisor
$K_{\pi}+a_1s_1+a_2s_2+\ldots +a_ns_n$ as any
section meeting the singularities has coefficient zero.  
We emphasize that the stability condition is the natural one arising
from the log minimal model program (cf. the proof of Proposition
\ref{prop:VCP}).  

The resulting moduli spaces $\ocM_{g,{\tilde \cA}}$ are easily described.
Let $\cA$ be the subsequence of ${\tilde \cA}$ containing all the
positive weights and assume that $|\cA|+N=|{\tilde \cA}|$.  
Each ${\tilde \cA}$-stable pointed curve consists of a $\cA$-stable
curve with $N$ additional arbitrary marked points, i.e., the points with
weight zero.  Hence we may identify
$$\ocM_{g,{\tilde \cA}}= \underbrace{\cC_{g,\cA} \times_{\ocM_{g,\cA}} \ldots 
\times_{\ocM_{g,\cA}} \cC_{g,\cA}}_{N \text{ times }},$$
so $\ocM_{g,{\tilde \cA}}$ is the $N$-fold fiber product of the universal
curve over $\ocM_{g,\cA}$.

The moduli spaces with zero weights differ from 
the original spaces in one crucial
respect:  they are generally singular.  For example, the local analytic equation
of a generic one-parameter deformation of a nodal curve may be written
$xy=t$, where $t$ is the coordinate on the base.  
The second fiber product of this family takes the form
$$x_1y_1=x_2y_2=t,$$
which is a threefold with ordinary double point. 

\subsubsection{Weights summing to two}
\label{subsect:wtstwo}
We restrict to the case $g=0$ and consider weight data
${\hat A}=(a_1,\ldots,a_n)$ where the weights are positive
rational numbers with $a_1+\ldots +a_n=2$.  Weighted pointed
curves of this type have previously been considered by Kawamata,
Keel, and McKernan in the context of the codimension-two subadjunction formulas 
(see \cite{Kaw2} and \cite{KeMc2}).  One can construct an explicit family
of such weighted curves over the moduli space
$$\cC({\hat A}) \ra \ocM_{0,n};$$
this family is realized as an explicit blow-down of the universal
curve over $\ocM_{0,n}$.  

In this paper we do not give a direct modular interpretation of spaces
$\ocM_{0,{\hat A}}$.  However, when each $a_j<1$ we may interpret
the geometric invariant theory quotient
$$({\bP}^1)^n//\SL_2$$       
with linearization $\cO(a_1,a_2,\ldots,a_n)$ as $\ocM_{0,{\hat A}}$  
(see Theorems \ref{thm:typical} and \ref{thm:atypical}).  
These spaces are often singular (see Remark \ref{rem:10sing}).      

\subsubsection{Weighted divisors}
We can also consider curves with weighted divisors rather than
weighted points.  A stable curve with
weighted divisors consists of a nodal connected curve $C$
of genus $g$,
a collection of effective divisors supported in its smooth locus
$$D_1,\ldots,D_m,$$
and positive weights $a_1,\ldots,a_m$, 
so that the sum $D:=a_1D_1+\ldots+a_mD_m$
has coefficient $\le 1$ at each point and $K_C+D$ is ample.
Writing $d_j=\deg(D_j)$, we can construct a moduli space
$$\oM_{g,((a_1,d_1),\ldots,(a_m,d_m))}$$
as follows.  We associate to this problem the weight data
$$\cB:=(\underbrace{a_1,\ldots,a_1}_{d_1\text{ times }},\ldots,
\underbrace{a_m,\ldots,a_m}_{d_m\text{ times}})$$
and the corresponding 
coarse moduli scheme of weighted pointed curves $\oM_{g,\cB}$.  
We take
$$\oM_{g,((a_1,d_1),\ldots,(a_m,d_m)}=
\oM_{g,\cB}/(S_{d_1}\times \ldots \times S_{d_m}),$$
where the product of symmetric groups acts componentwise
on the $m$ sets of sections.  We will not discuss the 
propriety of writing the moduli space as such a quotient,
except to refer the reader to chapter 1 of \cite{GIT}.

\section{Construction of the moduli space} \label{sect:construct}

\subsection{Preliminaries on linear series}
In this section we work over an algebraically closed
field $F$.  Given a curve $C$ and a smooth point $s$, note that
the ideal sheaf $I_s$ is invertible.  We write
$L(s)$ for $L\otimes I_s^{-1}$.  

\begin{prop}\label{prop:nosection}
Let $C$ be a connected nodal proper curve.
Let $M$ be an invertible sheaf such
that $M^{-1}$ is nef.  
Then $h^0(M)\neq 0$ if and only if $M$ is trivial.
\end{prop}
The nef assumption means $\deg(M|C_j)\le 0$
for each irreducible component $C_j\subset C$. 
  
\noindent {\em proof:} This is elementary if $C$ is smooth.  
For the general case, consider the normalization  
$\nu:C^{\nu}\ra C$, with irreducible components
$C^{\nu}_1,\ldots,C^{\nu}_N$.  We have the formula
$$p_a(C)=\sum_{j=1}^N p_a(C^{\nu}_j)+ \Delta - N +1$$
relating the arithmetic genera and $\Delta$,
the number of singularities of $C$.  
Recall the exact sequence
$$0\ra T \ra \Pic(C) \stackrel{\nu^*}{\ra} \Pic(C^{\nu})$$
where $T$ is a torus of rank $\Delta-N+1$.  
To reconstruct $M$ from $\nu^*M$, for each singular point $p\in C$
and points $p_1,p_2\in  C^{\nu}$ lying over $p$ we
specify an isomorphism $(\nu^*M)_{p_1}\simeq (\nu^*M)_{p_2}$,
unique up to scalar multiplication on $\nu^*M$.  In particular,
we obtain an exact sequence
$$0 \ra {\bG}_m \ra {\bG}_m^N \ra \bG_m^{\Delta} \ra T \ra 0.$$
If $M$ has a nontrivial section then 
$\nu^*M\simeq \cO_{C^{\nu}}$ and the section pulls back
to a section of $\cO_{C^{\nu}}$ constant and nonzero on each component.
Thus the corresponding element of $T$ is trivial and $M\simeq \cO_C$.
$\square$

\begin{prop}\label{prop:degrees}
Let $C$ be an irreducible nodal curve with arithmetic genus $g$, 
$B$ and $D$ effective
divisor of degrees $b$ and $d$ supported in the smooth locus of $C$.  
Let $M$ be an ample invertible sheaf that may be written
$$M=\omega_C^k(kB+D)\quad k>0,$$
and $\Sigma\subset C$ a subscheme of length $\sigma$ contained in
the smooth locus of $C$.  

Assume that $\sigma \le 2$ and $N\ge 4$ (resp. $\sigma\le 1$ and
$N\ge 3$).
Then $H^0(\omega_C(B+\Sigma)\otimes M^{-N})=0$.  This holds for $N=3$ 
(resp. $N=2$) except in the cases 
\begin{enumerate}
\item
$d=0,k=1,g=0,$ and $b=3$; or
\item
$d=0,k=1,g=1,b=1,$ and $\cO_C(\Sigma)\simeq \cO_C(\sigma B)$. 
\end{enumerate}
In these cases, all the sections are constant.
Finally, $H^0(\omega_C(B)\otimes M^{-N})=0$ when $N\ge 2$.          
\end{prop}

\noindent {\em proof:}  Setting $F=\omega_C(B+\Sigma)\otimes M^{-N}$, 
we compute
$$\deg(F)=(2g-2+b)+\sigma-N \deg(M)=(1-Nk)(2g-2+b)+\sigma-Nd.$$
We determine when these are nonnegative.  First assume $N\ge 3$
and $\sigma\le 2$.
Using the first expression for $\deg(F)$ and $\deg(M)\ge 1$, 
we obtain $2g-2+b\ge 1$.  The second expression implies that $2g-2+b\le 1$.
Thus $2g-2+b=1$, and the first expression gives
$N=3,\sigma=2,$ and $\deg(M)=1$;  the second expression yields
$k=1$ and $d=0$.  These are
the exceptional cases above.  

Now assume $N\ge 2$ and $\sigma\le 1$.  Repeating the
argument above, we find $2g-2+b=1$ and therefore $N=2,\sigma=1,
\deg(M)=1, k=1,$ and $d=0$.  Again, we are in one of the two
exceptional cases.  Finally, if $N\ge 2$ and $\sigma=0$, 
we obtain $(2g-2+b)>0$ from the first expression and $(2g-2+b)<0$ from
the second.  This proves the final assertion. 
$\square$

\begin{prop} \label{prop:series}
Let $C$ be a connected nodal curve of genus $g$, $D$ an effective divisor
supported in the smooth locus of $C$,
$L$ an invertible sheaf with $L\simeq\omega_C^k(D)$ for $k>0$.
\begin{enumerate}
\item
If $L$ is nef and $L\neq \omega_C$   
then $L$ has vanishing higher cohomology.
\item 
If $L$ is nef and has positive degree
then $L^N$ is basepoint free for $N\ge 2$.      
\item
If $L$ is ample then
$L^N$ is very ample when $N\ge 3$. 
\item
Assume $L$ is nef and has positive degree
and let $C'$ denote the image of $C$ under $L^N$ with $N\ge 3$.
Then $C'$ is a nodal curve with the same arithmetic genus as $C$,
obtained by collapsing the irreducible 
components of $C$ on which $L$
has degree zero.  Components on which $L$ has positive degree
are mapped birationally onto their images.     
\end{enumerate} 
\end{prop}
Our argument owes a debt to
Deligne and Mumford (\cite{DM} \S 1).

\noindent {\em proof:}  For the first statement, we use Serre duality
$h^1(L)=h^0(\omega_C\otimes L^{-1})$ and Proposition
\ref{prop:nosection} applied to $M=\omega_C\otimes L^{-1}$.  
One verifies easily that 
$$M^{-1}=L \otimes \omega_C^{-1} = \omega_C^{k-1}(D)$$ 
is the sum of a nef and an effective divisor.  
  
\

We prove the basepoint freeness statement.  Decompose
$$C=Z\cup_T C_+$$
where $Z$ contains the components on which $L$ has degree zero,
$C_+$ the components on which $L$ is ample, and $T$ is their
intersection.  Each connected component $Z_j \subset Z$ is a chain of
${\bP}^1$'s and has arithmetic genus zero.  A component $Z_j$ is
{\em type I} (resp. {\em type II}) if it contains one point
$t_j\in T$ (resp. two points $t_j',t_j''\in T$).  
    
It suffices to show that for each $p\in C$  
$$h^0(L^N\otimes I_p)=h^0(L^N)-1$$
where $N\ge 2$.  The vanishing assertion guarantees that
$L^N$ has no higher cohomology.  It suffices then to show that
$L^N\otimes I_p$ has no higher cohomology, or dually,
$\Hom(I_p,\omega_C\otimes L^{-N})=0$.  If $p$ is a smooth point then 
$I_p$ is locally free and
$$\Hom(I_p,\omega_C\otimes L^{-N})=H^0(\omega_C(p)\otimes L^{-N}).$$ 
We analyze the restriction of $F:=\omega_C(p)\otimes L^{-N}$ to the 
components of $C$.  We first restrict to $Z$:
$$h^0(Z_j,F|Z_j)=\begin{cases}
	0 \text{ if } p\not \in Z_j \text{ and } Z_j 
		\text{ is of type I;}\\
	1 \text{ if } p \in Z_j \text{ and }
	Z_j \text{ is of type I;}\\
	1 \text{ if } p \not \in Z_j \text{ and }
	Z_j \text{ is of type II;}\\
      2 \text{ if } p \in Z_j \text{ and }
	Z_j \text{ is of type II.}\\
\end{cases}$$
In each case
the sections are zero if they are zero at $T\cap Z_j$.
For components in $C_+$ we apply Proposition \ref{prop:degrees},  
where $M$ is restriction of $L$ to some irreducible component,
$\sigma=p$, and $B$ is the conductor.  (In what follows, on applying 
Proposition \ref{prop:degrees} we always assume $B$ contains
the conductor.)  The Proposition gives that the restriction to each 
component has no nontrivial sections, except perhaps when $p$ lies
on a component $E$ listed in the exceptional cases.  Then 
the sections of $F|E$ are constant and $E\simeq {\bP}^1$
because $p\not \in B$.  Thus provided
$p$ is not contained in a component $E\simeq \bP^1 \subset C_+$ with
$|B|=3$, we obtain that $h^0(C,F)=0$.  Indeed, clearly
$h^0(F|C_+)=0$ and the analysis of cases above 
yields $h^0(C,F)=0$.  If $p$ does sit on such 
a component $E$, then $C$ contains a component of type I or
a second irreducible component of $C_+$;  the conductor $B\subset E$ 
has three elements, so there must be at least one other non-type II
component.  The restriction of $F$ to such a component has only
trivial sections, and since the restrictions to all the
other components have at most constant sections,
we conclude $h^0(C,F)=0$.   

If $p$ is singular, let $\beta:\hC \ra C$ be the
blow-up of $C$ at $p$ and $p_1,p_2\in \hC$ the points
lying over $p$, so that
$$\Hom(I_p,\omega_C\otimes L^{-N})=H^0(\beta^*(\omega_C\otimes L^{-N}))=
H^0(\omega_{\hC}(p_1+p_2)\beta^*L^{-N}).$$
We write $\hF=\omega_{\hC}(p_1+p_2)\beta^*L^{-N}$.
As before, we decompose
\begin{equation}
\hC={\hC}_+ \cup_{\hT} \hZ \label{eqn:decomp}
\end{equation} 
where $L$ is positive on $\hC_+$, and use
$\hZ_j$ to denote a connected component of $Z$.  Note that $\beta(\hZ)=Z$
and the $p_i$ are not both contained on some $\hZ_j$;  otherwise,
$\beta(\hZ_j)$ would have positive arithmetic genus.  Similarly, neither
of the $p_i$ lie on a type II component (with respect to
decomposition (\ref{eqn:decomp})).  It follows that $h^0(\hF|\hZ_j)$ has dimension
at most one, and any section vanishes if it vanishes along
$\hT\cap \hZ_j$.  For components  
of $\hC_+$, we apply Proposition \ref{prop:degrees}, with $p_1,p_2\in B$
and $\sigma=0$, to show that $h^0(\hF|\hC_+)=0$.  Again, we conclude
that $h^0(\hC,\hF)=0$.      

\

For ampleness, take $p$ and $q$ to be points of $C$, 
not necessarily smooth or distinct, with
ideal sheaves $I_p$ and $I_q$.  Again, it suffices to
prove $h^0(I_pI_qL^N)=0$ for $N\ge 3$, or dually,
$\Hom(I_pI_q,\omega_C\otimes L^{-N})=0$.
When $p$ and $q$ are smooth points the assertion follows as above
from Proposition \ref{prop:degrees}, again with $B$ as the conductor.    
The only case requiring additional argument is when $p$ and $q$ both lie
on a component corresponding to one of the exceptional cases.  Again,
the nontrivial sections on this component are constants, whereas we have
only zero sections on the other components.  

We may assume that $p$
is singular and write $\beta:\hC\ra C$ for the blow-up at $p$.  Then
for each invertible sheaf $R$ on $C$ we have
$$\Hom(I_p^2,R)=H^0(\beta^*R(p_1+p_2))$$
where $p_1$ and $p_2$ are the points of $\hC$ lying over $p$.      
When $q$ is smooth we obtain
$$\Hom(I_pI_q,\omega_C\otimes L^{-N})
=H^0(\hC,\omega_{\hC}(p_1+p_2+q) \otimes \beta^*L^{-N}).$$
Note that the restriction of $\beta^*L$ to each component
still satisfies the hypothesis of Proposition
\ref{prop:degrees};  we take $p_1,p_2\in B$ and $q=\Sigma$.  
When $q$ is singular and disjoint from $p$ we obtain        
$$\Hom(I_pI_q,\omega_C\otimes L^{-N})
=H^0(\hC,\omega_{\hC}(p_1+p_2+q_1+q_2) \otimes \beta^*L^{-N}),$$
where $\hC \ra C$ is the blow-up at $p$ and $q$.
We apply Proposition \ref{prop:degrees} with $p_1,p_2,q_1,q_2\in B$
and $\Sigma=\emptyset$.  When $p=q$ we obtain    
$$\Hom(I_p^2,\omega_C\otimes L^{-N})
=H^0(\hC,\omega_{\hC}(2p_1+2p_2) \otimes \beta^*L^{-N}),$$
and we apply Proposition \ref{prop:degrees}
with $p_1,p_2\in B$ and $\Sigma=\{p_1,p_2\}$, 
accounting for the exceptional cases as before.  

\

We prove the last assertion.  We have an exact sequence
$$0\ra H^0(C,L^N) \ra H^0(C_+,L^N|C_+) \oplus H^0(Z,L^N|Z) \ra H^0(T,L^N|T).$$  
Choose $N$ large so that $(L^N|C_+)(-T)$ is very ample.  
Since $L^N|Z$ is trivial, the image of $C$ under $L^N$ is
obtained from $C_+$ by identifying pairs of points in $T$
corresponding to $t_j',t_j''$ in some $Z_j$,i.e.,
by `collapsing'
each $Z_j$ to a point.  Let $C'$ denote the resulting curve, which has
the same arithmetic genus as $C$, and $r:C \ra C'$ the resulting map.  
Writing $D'=r(D)$, we have $r^*\omega^k_{C'}(D')=L$.  Thus the 
sections of $L^N$ ($N\ge 3$) induce $r$, as $(\omega^k_{C'}(D'))^N$
is very ample on $C'$.  
$\square$

\begin{rem} \label{rem:typesnew}
The adjunction formula gives
precise information about the points of $D$ lying in $Z$.
Each connected component $Z_j$ of type II is disjoint from $D$.
Recall that components $Z_j$ of type I are chains 
$Z_{j1}\cup \ldots \cup Z_{jm}$ of ${\bP}^1$'s,
intersecting $C_+$ in a point of one of the ends of the
chain (say $Z_{j1}$).  Then $D|Z_j$ is supported in 
the irreducible component $Z_{jm}$ at the opposite end
and $\deg(D|Z_{jm})=k$.  
\end{rem}     

Applying Proposition \ref{prop:series} and Remark 
\ref{rem:typesnew} with $D=k(b_1s_1+\ldots+b_ns_n)$,
we obtain the following:
\begin{cor} \label{cor:section}
Let $(C,s_1,\ldots,s_n)$ be a nodal pointed curve of genus $g$,  
$b_1,\ldots,b_n$ nonnegative rational numbers, and $k$ a 
positive integer such that each $kb_i$ is integral.  
Assume that $L:=\omega_C^k(b_1s_1+\ldots+b_ns_n)$ is nef and has
positive degree.  

For $N\ge 3$ the sections of $L^N$ induce
a dominant morphism $r:C \ra C'$ to a nodal curve of genus $g$.
This morphism collapses irreducible components of $C$ on which
$L$ has degree zero, and maps the remaining components birationally
onto their images. 
If $\cB=(b_{i_1},\ldots,b_{i_r}), i_1<i_2< \ldots <i_r$,
denotes the set of all nonzero weights and $s'_i=r(s_i)$, 
then $(C',s'_{i_1},\ldots,s'_{i_r})$ is a stable pointed curve
of type $(\cB,g)$.  
\end{cor}  

\noindent {\em proof: } The only claim left to verify is
the singularity condition.  Each $s_i$ lying in a component of
type II necessarily has weight $b_i=0$.  
Thus no points with positive weight
are mapped to singularities of the image $C'$. 
The points $\{s_{j_1},\ldots,s_{j_a}\}$ 
lying on a single component of type I have weights summing to one, 
i.e., $b_{j_1}+\ldots +b_{j_a}=1$.  $\square$   

We also obtain the following relative statement:
\begin{thm}\label{thm:logcan}
Let $\pi:(C,s_1,\ldots,s_n)\ra B$ be a family of nodal pointed curves
of genus $g$,
$b_1,\ldots,b_n$ nonnegative rational numbers, and $k$ a positive
integer such that each $kb_i$ is integral.
Set $L=\omega^k_{\pi}(kb_1s_1+\ldots +kb_ns_n)$ and assume that $L$ is
$\pi$-nef and has positive degree.  
For $N\ge 3$, 
$\Proj(\oplus_{m\ge 0}\pi_*L^{mN})$ defines a flat family of nodal
curves $C'$ with sections $s'_1,\ldots,s'_n$.  
If $\cB=(b_{i_1},\ldots,b_{i_r}), i_1<i_2< \ldots <i_r,$
denotes the set of all nonzero weights, 
then $(C',s'_{i_1},\ldots,s'_{i_r})$ is a family of stable pointed curves
of type $(\cB,g)$.       
\end{thm}       
The new family may be considered as the {\em log canonical model} of $C$
relative to $K_{\pi}+a_1s_1+\ldots +a_ns_n$.  

\noindent {\em proof:}  The vanishing assertion of Proposition \ref{prop:series}
implies the formation of $\pi_*L^N$ commutes with base extensions
$B'\ra B$.  Hence we may apply the fiberwise assertions
of Corollary \ref{cor:section}
to $(C',s'_1,\ldots,s'_n)\ra B$.  We therefore
obtain a family of pointed stable curves of type $(\cB,g)$.
$\square$

\subsection{The log minimal model program and the valuative
criterion}
To prove that our moduli problem is proper, we shall apply the valuative
criterion for properness (cf. \cite{LM} 7.5).  
The most important step is the following:
\begin{prop}\label{prop:VCP}
Let $R$ be a DVR with quotient field $K$, $\Delta=\Spec R$, $\Delta^*=\Spec K,$ 
$(g,\cA)$ a collection of input data, 
$\pi^*:(C^*,s^*_1,\ldots,s^*_n)\ra \Delta^*$ a family of stable pointed curves of
type $(g,\cA)$.  Then there exists the spectrum of a DVR ${\tilde \Delta}$,
a finite ramified morphism ${\tilde \Delta} \ra \Delta$, and a family
$\pi^c: (\cC^c,s_1^c,\ldots,s^c_n)\ra {\tilde \Delta}$ of stable pointed curves of 
type $(g,\cA)$, isomorphic to 
$$(C^*\times_{\Delta}{\tilde \Delta},s_1^*\times_{\Delta}{\tilde \Delta},\ldots,
s_n^*\times_{\Delta}{\tilde \Delta})$$
over ${\tilde \Delta}^*$.  The family $(\cC^c,s_1^c,\ldots,s_n^c)$ is unique
with these properties.
\end{prop}
{\em proof:}  We first reduce to the case where $C^*$ is geometrically normal
with disjoint sections.  If sections $s_{i_1},\ldots,s_{i_r}$ coincide over
the generic point, we replace these by a single section with weight
$a_{i_1}+\ldots+a_{i_r}$.    
Choose a finite extension of $K$ over which each irreducible component of $C^*$
is defined, as well as each singular point.  Let $C^*({\nu})$ be the normalization, 
$s^*({\nu})_1,\ldots,s^*({\nu})_n$ the proper transforms of the sections, and
$s^*({\nu})_{n+1},\ldots,s^*({\nu})_{n+b}$ the points of the conductor.  Then
$(C^*({\nu}),s^*({\nu})_1,\ldots,s^*({\nu})_{n+b})$ is stable with respect to the weights
$(\cA,1,\ldots,1)$.  Once we have the stable reduction of $C^*(\nu)$, the stable
reduction of $C^*$ is obtained by identifying corresponding pairs of points of the
conductor.
    
Applying the valuative criterion for properness for $\cM_{g,n}$
(which might entail a base-change ${\tilde \Delta} \ra \Delta$), 
we reduce to the case where $(C^*,s_1^*,\ldots,s^*_n)$ extends to a family 
$\pi:(C,s_1,\ldots,s_n)\ra \Delta$
of stable curves in $\cM_{g,n}$.  If this family is stable with respect to 
the weight data $\cA$ (i.e., $K_{\pi}+a_1s_1+\ldots a_ns_n$ is ample
relative to $\pi$) then there is nothing to prove.  
We therefore assume this is not the case.

Our argument uses the log minimal model program
to obtain a model on which our log canonical divisor is ample.   
This is well-known for surfaces over fields of arbitrary characteristic
(see Theorem 1.4 of \cite{Fu} or \cite{KK}), but perhaps
less well-known in the mixed characteristic case.  For
completeness, we sketch a proof.

Let $\lambda$ be the largest number for which
$$D:=\lambda (K_{\pi}+s_1+\ldots+s_n) + (1-\lambda)(K_{\pi}+a_1s_1+\ldots
+a_ns_n)$$
{\em fails} to be ample.  Our assumptions imply $0\le \lambda <1$; 
$\lambda$ is rational because there are only finitely many 
(integral projective) curves lying in fibers of $\pi$.  
Our argument is by induction on the number of such curves.
The ${\bQ}$-divisor $D$ is nef and has positive degree, and we choose  
$L$ to be the locally free sheaf associated to a suitable
multiple of $D$.  Applying Theorem \ref{thm:logcan},
we obtain a new family of pointed curves 
$\pi':(C',s_1'+\ldots+s'_n)\ra {\tilde \Delta}$,
agreeing with the original family away from the central fiber, and
stable with respect to the weight data 
$\cB(1):=\lambda(1,1,\ldots,1)+(1-\lambda)\cA$.  Note that
in passing from $C$ to $C'$, we have necessarily contracted
some curve in the central fiber $\pi$.  Now if $\cB(1)=\cA$ 
(i.e., if $\lambda=0$) the proof is complete.  Otherwise,
we repeat the procedure above using a suitable log divisor 
$$D':=\lambda (K_{\pi'}+b_1(1)s'_1+\ldots+b_n(1)s'_n) + (1-\lambda)(K_{\pi'}+a_1s'_1+\ldots
+a_ns'_n).$$
We continue in this way until $\cB(j)=\cA$;  this process
terminates because there are only finitely many curves
in the central fiber to contract.  $\square$
  
\subsection{Deformation theory}
Let $(C,s_1,\ldots,s_n)$ be a weighted pointed stable curve of genus $g$
with weight data $\cA$, defined over a field $F$.  
We compute its infinitesimal automorphisms and 
deformations.  We regard the pointed curve as a map 
$$s:\Sigma \ra C,$$
where $\Sigma$ consists of $n$ points, each mapped to the corresponding $s_j\in C$.

The infinitesimal deformation theory of maps 
was analyzed by Ziv Ran \cite{Ra}; he worked with 
holomorphic maps of reduced analytic spaces, 
but his approach also applies in an algebraic context.
We recall the general formalism.  
Infinitesimal automorphisms, deformations, and obstructions
of $s:\Sigma\ra C$ are denoted 
by $T^0_s,T^1_s,$ and $T^2_s$ respectively.
Similarly, we use $T^i_C=\Ext^i_C(\Omega^1_C,\cO_C)$ 
and $T^i_{\Sigma}=\Ext^i_{\Sigma}(\Omega^1_{\Sigma},\cO_{\Sigma})$ 
for the analogous groups associated to $C$ and $\Sigma$.
Finally, we consider the mixed group $$\Hom_s(\Omega^1_C,\cO_{\Sigma})=\Hom_C(\Omega^1_C,s_*\cO_{\Sigma})=
\Hom_{\Sigma}(s^*\Omega^1_C,\cO_{\Sigma})$$
and the associated $\Ext$-groups,
denoted $\Ext^i_s(\Omega^1_C,\cO_{\Sigma})$ and computed by
either of the spectral sequences
$$E^{p,q}_2=\Ext^p_C(\Omega^1_C,\bR^qs_* \cO_{\Sigma}) \quad
E^{p,q}_2=\Ext^p_{\Sigma}(\bL^qs^* \Omega^1_C,\cO_{\Sigma}).$$
We obtain long exact sequences
\begin{eqnarray*}
0\ra T^0_s \ra T^0_C\oplus T^0_{\Sigma} \ra \Hom_s(\Omega^1_C,\cO_{\Sigma}) 
 \ra T^1_s \ra T^1_C \oplus T^1_{\Sigma}\\
\ra \Ext^1_s(\Omega^1_C,\cO_{\Sigma}) 
 \ra T^2_s \ra T^2_C \oplus T^2_{\Sigma}
\ra \Ext^2_s(\Omega^1_C,\cO_{\Sigma}).
\end{eqnarray*}

In our situation $T^i_{\Sigma}=0$ 
($\Sigma$ is reduced zero-dimensional), $T^2_C=0$ 
($C$ is a nodal curve), $\Ext^i_s(\Omega_C,\cO_{\Sigma})=
\Ext^i_{\Sigma}(s^*\Omega_C,\cO_{\Sigma})$ ($\Omega^1_C$
is free along $s(\Sigma)$), and thus $\Ext^i_s(\Omega_C,\cO_{\Sigma})=0$
for $i>0$.  Hence the exact sequence boils down to
$$0\ra T^0_s \ra \Hom_C(\Omega^1_C,\cO_C) \ra 
\Hom_{\Sigma}(s^*\Omega^1_C,\cO_{\Sigma}) \\
\ra T^1_s \ra \Ext^1_C(\Omega^1_C,\cO_C) \ra 0$$
and $T^2_s=0$.    

Let $D\subset C$ denote the support of $a_1s_1+\ldots+a_ns_n$.
Note that 
$$D-(a_1s_1+\ldots+a_ns_n)$$ 
is an effective ${\bQ}$-divisor and so the positivity condition
guarantees that $\omega_C(D)$ is ample. 
The map $s$ factors
$$s:\Sigma \ra D \hookrightarrow C$$
which gives a factorization
$$\Hom_C(\Omega^1_C,\cO_C) \ra \Hom_C(\Omega^1_C,\cO_D)
\ra \Hom_{\Sigma}(s^*\Omega^1_C,\cO_{\Sigma}).$$
The second step is clearly injective.  The kernel of the
first step is 
$$\Hom_C(\Omega^1_C,\cO_C(-D))\simeq H^0((\omega_C(D))^{-1})=0.$$  
Thus $T^0_s=0$ and $T^1_s$ has dimension $3g-3+n$.  
We summarize this in the following proposition:
\begin{prop}\label{prop:infdef}
Let $(C,s_1,\ldots,s_n)$ be a weighted pointed stable curve of genus $g$
with weight data $\cA$.  Then this curve admits no infinitesimal
automorphisms and its infinitesimal deformation space is unobstructed of
dimension $3g-3+n$.
\end{prop}

\subsubsection{The canonical class}

We digress to point out consequences of this analysis
for the moduli stack.
The tangent space to $\ocM_{g,\cA}$ at 
$(C,s_1,\ldots,s_n)$ sits in the exact sequence
\begin{eqnarray*}
0\ra \Hom_C(\Omega^1_C,\cO_C) \ra 
\oplus_{j=1}^n \Hom_C(\Omega^1_C,\cO_{s_j}) \\
\ra  T_{(C,s_1,\ldots,s_n)}\ocM_{g,\cA} 
\ra \Ext_C^1(\Omega^1_C,\cO_C) \ra 0.
\end{eqnarray*}
The cotangent space sits in the dual exact sequence
$$0\ra H^0(\Omega^1_C\otimes \omega_C) \ra 
T^*_{(C,s_1,\ldots,s_n)}\ocM_{g,\cA} \ra 
\oplus_{j=1}^n \Omega^1_C|s_j \ra H^1(\Omega^1_C\otimes \omega_C)
\ra 0.$$
Now let
$$\pi:\cC_{g,\cA}\ra \ocM_{g,\cA}$$
be the universal curve and $s_j$ the corresponding sections. 
The exact sequences above cannot be interpretted as exact sequences
of vector bundles on the moduli stack because  
$h^0(\Omega^1_C\otimes \omega_C)$ and $h^1(\Omega^1_C\otimes \omega_C)$
are nonconstant.  However, there is a relation in the derived category
which, on combination with the Grothendieck-Riemann-Roch Theorem, yields a formula
for the canonical class of the moduli stack (cf. \cite{HM} pp. 159). 
This takes the form
$$K_{\ocM_{g,\cA}}=\frac{13}{12}\kappa'(\cA) - \frac{11}{12}\nu(\cA) + 
\sum_{j=1}^n \psi_j(\cA),$$
where $\kappa'(\cA)=\pi_*[c_1(\omega_{\pi})^2]$,
$\nu(\cA)$ the divisor parametrizing nodal curves,
and $\psi_j(\cA)=c_1(s_j^*\omega_{\pi})$.  

\subsubsection{An alternate formulation}
We sketch an alternate formalism for the deformation theory
of pointed stable curves.  This was 
developed by Kawamata \cite{Kaw1} in an analytic context
for logarithmic pairs $(X,D)$ consisting of a proper nonsingular 
variety $X$ and a normal crossings divisor $D\subset X$.  
This approach is more appropriate
when we regard the boundary as a {\em divisor} rather than
the union of a sequence of sections.  In particular, this
approach should be useful for higher-dimensional generalizations
of weighted pointed stable curves, like stable log surfaces.

We work over an algebraically closed field $F$.
Let $X$ be a scheme and $D_1,\ldots,D_n$ a sequence
of distinct effective Cartier divisors (playing the role of
the irreducible components of the normal crossings
divisor).  We define the sheaf $\Omega^1_X\left<D_1,\ldots,D_n\right>$
of differentials on $X$ 
with logarithmic poles along the collection $D_1,\ldots,D_n$.  
Choose an open affine subset
$U\subset X$ so that each $D_j$ is defined by an equation
$f_j\in \cO_U$.  Consider the module           
$$\Omega^1_U\left<D_1,\ldots,D_n\right>:=
(\Omega^1_U\oplus \Omega^1_U e_{f_1}\oplus \ldots \oplus \Omega^1_U e_{f_n})/
\left<df_j-f_j e_{f_j}\right>.$$
Up to isomorphism, this is independent of the choice of the $f_j$;
indeed, if $f_j=ug_j$ with 
$u\in \cO_U^*$ then we have the substitution $e_{g_j}=e_{f_j}-du/u$.     
The inclusion by the first factor induces a natural injection
$$\Omega^1_U \hookrightarrow \Omega^1_U\left<D_1,\ldots,D_n\right>$$
with cokernel $\oplus_{j=1}^n\cO_{D_j}$, where each summand is generated
by the corresponding $e_{f_j}$.  We therefore obtain the following 
natural exact sequence of $\cO_X$-modules
$$0 \ra \Omega^1_X \ra \Omega^1_X\left<D_1,\ldots,D_n\right> \ra 
\oplus_{j=1}^n\cO_{D_j}\ra 0.$$

In the case where $X$ is smooth and the $D_j$ are smooth, reduced, 
and meet in normal crossings, 
we recover the standard definition of differentials with 
logarithmic poles and the exact sequence is the ordinary residuation
exact sequence.  When $D_j$ has multiplicity $>1$, the sheaf 
$\Omega^1_X\left<D_1,\ldots,D_n\right>$ has torsion along $D_j$.        
Note that the residuation exact sequence is split when the multiplicities
are divisible by the characteristic.

Assume for simplicity that $X$ is projective and 
smooth along the support of the $D_j$.  We claim
that first-order deformations of $(X,D_1,\ldots,D_n)$ correspond
to elements of
$$\Ext^1_X(\Omega^1_X\left<D_1,\ldots,D_n\right>,\cO_X).$$
The resolution
$$0\ra \cO_X(-D_j) \ra \cO_X \ra \cO_{D_j} \ra 0$$ 
implies $\Ext^1_X(\cO_{D_j},\cO_X)=H^0(\cO_{D_j}(D_j)),$
so the residuation exact sequence yields
$$\oplus_{j=1}^n H^0(\cO_{D_j}(D_j)) \ra 
\Ext^1_X(\Omega^1_X\left<D_1,\ldots,D_n\right>,\cO_X) \ra 
\Ext^1_X(\Omega^1_X,\cO_X).$$
Of course, $H^0(\cO_{D_j}(D_j))$ parametrizes first-order deformations of $D_j$
in $X$ and $\Ext^1_X(\Omega^1_X,\cO_X)$ parametrizes first-order
deformations of the ambient variety $X$.

\subsection{Construction of the moduli stack}
In this section, we prove all the assertions of
Theorem \ref{thm:existence} except the projectivity
of the coarse moduli space, which will be proved
in \S \ref{subsect:proj}.  For a good general 
discussion of how moduli spaces are constructed,
we refer the reader to \cite {DM} or \cite{Vi}.

We refer to \cite{LM}, \S 3.1 and 4.1, for the definition of
a stack.  The existence of the moduli space as a stack
follows from standard properties of descent:
families of stable pointed curves of type
$(g,\cA)$ satisfy effective descent and
$\Isom$ is a sheaf in the \'etale topology.

We introduce an `exhausting
family' for our moduli problem, i.e., a scheme which is an atlas
for our stack in the smooth topology.
Set 
$$L=\omega_{\pi}^k(ka_1s_1+\ldots + ka_ns_n) $$
where $k$ is the smallest positive integer such that each
$ka_j$ is integral.    
Let $d=\deg(L^3)=3k(2g-2+a_1+\ldots+a_n)$
and consider the scheme $H_0$ 
parametrizing $n$-tuples $(s_1,\ldots,s_n)$ in $\bP^d$, and 
the scheme $H_1$ parametrizing
genus $g$, degree $d$ curves $C\subset {\bP}^d$.  Let
$$U\subset H_0 \times H_1$$
be the locally closed subscheme satisfying the following
conditions:
\begin{enumerate}
\item{$C$ is reduced and nodal;}
\item{$s_1,\ldots,s_n \subset C$ and is contained in the smooth locus;}
\item{$\cO_C(+1)=L^3$.}
\end{enumerate}
We are using the fact that two line bundles
(i.e., $\cO_C(+1)$ and $L^{\otimes 3}$) coincide on a locally 
closed subset.

We shall now prove that our moduli stack is algebraic using
Artin's criterion (see \cite{LM} \S 10.) 
Proposition \ref{prop:series} implies that each curve
in $\ocM_{g,\cA}$ is represented in $U$.  Furthermore,
isomorphisms between pointed curves
in $U$ are induced by projective equivalences in $\bP^d$.
It follows that $U\ra \ocM_{g,\cA}$ is smooth and surjective,
and that the diagonal 
$\ocM_{g,\cA}\ra \ocM_{g,\cA}\times \ocM_{g,\cA}$ is representable,
quasi-compact, and separated.  We conclude that 
$\ocM_{g,\cA}$ is an algebraic stack of finite type.

To show that the stack is 
Deligne-Mumford, it suffices to show that our pointed
curves have `no infinitesimal automorphisms,'
i.e., that the diagonal
is unramified \cite{LM}, 8.1.  
This follows from Proposition \ref{prop:infdef}.
The moduli stack is proper over $\bZ$ by the valuative
criterion for properness (cf. \cite{LM} 7.5) and Proposition \ref{prop:VCP}.
The moduli stack is smooth over $\bZ$ if, for each curve defined 
over a field, the 
infinitesimal deformation space is unobstructed. 
This also follows from Proposition \ref{prop:infdef}.  
$\square$

\subsection{Existence of polarizations}
\label{subsect:proj}

We now construct polarizations for the moduli spaces of weighted pointed 
stable curves, following methods of Koll\'ar \cite{Ko1} (see also \cite{KoMc}).  
We work over an algebraically closed field
$F$.  The first key concept is the notion of a semipositive sheaf.  
Given a scheme (or algebraic
space) $X$ and a vector bundle $E$ on $X$,
we say $E$ is {\em semipositive} if
for each complete curve $C$ and map $f:C\ra X$,
any quotient bundle of $f^*E$ has nonnegative degree.  
Second, we formulate precisely what it means to say
that the `classifying map is finite.'
Given an algebraic group $G$, a $G$-vector bundle $W$ on $X$ of rank $w$
and a quotient vector bundle $Q$ of rank $q$, the classifying map
should take the form
$$u:X \ra \Gr(w,q)/G$$
where the Grassmannian denotes the $q$-dimensional quotients of 
fixed $w$-dimensional space.  Since the orbit space need
not exist as a scheme, we regard $u$ as a set-theoretic map
on closed points $X(F)\rightarrow \Gr(w,q)(F)/G(F)$.    
We say that the classifying map $u$ is finite when it has finite fibers
and each point of the image has finite stabilizer.   

The following result, a slight modification of the 
Ampleness Lemma of \cite{Ko1}, 
allows us to use semipositive sheaves to construct polarizations:
\begin{prop}
\label{prop:amplem}
Let $X$ be a proper algebraic space, $W$ a semipositive vector bundle
with structure group $G$ and rank $w$.  
Let $Q_1,\ldots,Q_m$ be quotient
vector bundles of $W$ with ranks $q_1,\ldots,q_m$.  Assume that
\begin{enumerate} 
\item{$W=\Sym^d(V)$ for some vector bundle $V$ of rank $v$
and $G=GL_v$;}
\item{the classifying map
$$u:X\ra (\prod_{j=1}^m \Gr(w,q_j))/G$$
is finite.}
\end{enumerate}
Then for any positive integers $c_1,\ldots,c_m$ the line bundle 
$$\det(Q_1)^{c_1}\otimes \det(Q_2)^{c_2} \ldots \otimes \det(Q_m)^{c_m}$$
is ample.
\end{prop}
In the original result $m=1$, so the classifying map
takes values in the quotient of a single Grassmannian rather than
a product of Grassmannians.  However, the argument of \cite{Ko1}
generalizes easily to our situation, so we refer to this paper for 
the details.  

To apply this result we need to produce semipositive vector 
bundles on the moduli space $\ocM_{g,\cA}$.  Choose 
an integer $k\ge 2$ so that each $ka_i$ is an integer.
Given a family $\pi:(C,s_1,\ldots,s_n)\ra B$ we take
$$V:=\pi_*[\omega_{\pi}^k(ka_1s_1+\ldots+ka_ns_n)],$$
which is locally free.  Indeed, a degree computation
(using the assumption $k\ge 2$) yields
$$\omega_{C_b}^k(ka_1s_1(b)+\ldots+ka_ns_n(b))\neq \omega_{C_b},$$
so Proposition \ref{prop:series} guarantees the vanishing of
higher cohomology.  The sheaf $V$ 
is semipositive on $B$ by Proposition 4.7 of \cite{Ko1}.    
It follows that each symmetric product $\Sym^d(V)$ is also
semipositive (\cite{Ko1} 3.2).

We may choose $k$ uniformly large 
so that $\omega_{\pi}^k(ka_1s_1+\ldots +ka_ns_n)$ 
is very ample relative to $\pi$ for {\em any} family $\pi$
(see Proposition \ref{prop:series}). 
We shall consider the multiplication maps
$$\mu_d:\Sym^d(V) \ra \pi_*[\omega_{\pi}^{dk}(dka_1s_1+\ldots +dka_ns_n)]$$    
and the induced restrictions
$$\Sym^d(V) \ra Q_j, \quad 
Q_j:=s^*_j[\omega_{\pi}^{dk}(dka_1s_1+\ldots +dka_ns_n)].$$
These are necessarily surjective and each has kernel consisting
of the polynomials vanishing at the corresponding section.    
We also choose $d$ uniformly large so that 
$$W:=\Sym^d(V) \ra Q_{n+1}, \quad
Q_{n+1}:=\pi_*[\omega_{\pi}^{kd}(kda_1s_1 + \ldots + kda_ns_n)],$$
is surjective and the fibers of $\pi$ are cut out by elements
of the kernel (i.e., they are defined by equations of degree $d$).
In particular, the pointed curve $(C_b,s_1(b),\ldots,s_n(b))$ can be recovered
from the image of $b$ under the classifying map associated to $W$ and its
quotients $Q_1,\ldots,Q_{n+1}$.  The stabilizer of the image 
corresponds to the automorphisms of 
$${\bP}(H^0(\omega_{C_b}^k(ka_1s_1(b)+\ldots+ka_ns_n(b))))$$
preserving the equations of $C_b$ and the sections.  This is 
finite since our curve has no infinitesimal
automorphisms (Proposition \ref{prop:infdef}).    

If $\oM_{g,\cA}$ admits a universal family the existence of a polarization follows from Proposition \ref{prop:amplem}.
In general, we obtain only a line bundle $\cL$ on the moduli stack 
$\ocM_{g,\cA}$, but some positive power $\cL^N$
descends to the scheme $\oM_{g,\cA}$.  This power is ample.    
Indeed, there exists a family 
$\pi:(C,s_1,\ldots,s_n) \ra B$ of curves in $\ocM_{g,\cA}$ so
that the induced moduli map $B\ra \oM_{g,\cA}$ is finite
surjective (cf. \cite{Ko1} \S 2.6-2.8).  
The bundle $\cL^N$ is functorial in the sense that
it pulls back to the corresponding product          
$$\det(Q_1)^{Nc_1}\otimes \det(Q_2)^{Nc_2} \ldots \otimes \det(Q_{n+1})^{Nc_{n+1}}$$
associated with our family.  This is
ample by Proposition \ref{prop:amplem}.  Our proof 
of Theorem \ref{thm:existence} is complete.$\square$

\section{Natural transformations}
\label{sect:reduction}

\subsection{Reduction and forgetting morphisms}
\begin{thm}[Reduction]\label{thm:reduction}
Fix $g$ and $n$ and let $\cA=(a_1,\ldots,a_n)$ and $\cB=(b_1,\ldots,b_n)$
be collections of weight data so that $b_j\le a_j$ for each $j=1,\ldots,n$.
Then there exists a natural birational reduction morphism
$$\rho_{\cB,\cA}:{\ocM}_{g,\cA} \ra
{\ocM}_{g,\cB}.$$
Given an element $(C,s_1,\ldots,s_n) \in \ocM_{g,\cA}$,
$\rho_{\cB,\cA}(C,s_1,\ldots,s_n)$ is obtained by successively
collapsing components of $C$
along which 
$K_C+b_1s_1+\ldots+b_ns_n$ fails to be ample.
\end{thm}
\begin{rem}
The proof of Theorem \ref{thm:reduction} also applies 
when some of the weights
of $\cB$ are zero (see \S \ref{subsect:zerowts}).  
\end{rem}
\begin{thm}[Forgetting]
Fix $g$ and let $\cA$ be a collection of weight data and
$\cA':=\{a_{i_1},\ldots,a_{i_r} \} \subset \cA$ a subset 
so that $2g-2+a_{i_1}+\ldots+a_{i_r}>0$.  
Then there exists a natural forgetting morphism
$$\phi_{\cA,\cA'}:{\ocM}_{g,\cA} \ra
{\ocM}_{g,\cA'}.$$
Given an element $(C,s_1,\ldots,s_n) \in \ocM_{g,\cA}$,
$\phi_{\cA,\cA'}(C,s_1,\ldots,s_n)$ is obtained by successively
collapsing components of $C$
along which 
$K_C+a_{i_1}s_{i_1}+\ldots+a_{i_r}s_{i_r}$ fails to be ample.
\end{thm}
We refer the reader to Knudsen and Mumford \cite{KnMu}\cite{Kn}\cite{Kn2}
for the original results on the moduli space of unweighted
pointed stable curves.

\noindent {\em proof:}  We shall prove these theorems simultaneously,
using $\psi_{\cA,\cB}$ to denote either a reduction or a
forgetting map, depending on the context.  
Let $(g,\cA)$ be a collection of input data, 
${\hat \cB}=(b_1,\ldots,b_n)\in \bQ^n$ so that 
$0\le b_j \le a_j$ for each $j$ 
and $2g-2+b_1+\ldots+b_n>0$.  Let 
$\cB=(b_{i_1},\ldots,b_{i_r}), i_1<i_2<\ldots<i_r,$ be obtained
by removing the entries of ${\hat \cB}$ which are zero.  
We shall define a natural transformation of functors
$$\psi_{\cB,\cA}:\ocM_{g,\cA}\ra \ocM_{g,\cB}.$$
Our construction will yield a morphism in the
category of stacks, which therefore induces a morphism
of the underlying coarse moduli schemes.  

Consider ${\hat \cB}(\lambda)=\lambda \cA+(1-\lambda){\hat \cB}$
for $\lambda \in \bQ$, $0<\lambda <1$, and write
${\hat \cB}(\lambda)=(b_1(\lambda),\ldots,b_n(\lambda))$.  
We may assume there
exists no subset $\{i_1,\ldots,i_r \} \subset \{1,\ldots, n \}$
such that $b_{i_1}(\lambda)+\ldots +b_{i_r}(\lambda)=1$.  
If $\cA$ and $\cB$ do not satisfy this assumption, then
there is a finite sequence
$1>\lambda_0>\lambda_1>\ldots >\lambda_N=0$
so that each $({\hat \cB}(\lambda_j),{\hat \cB}(\lambda_{j+1}))$
does satisfy our assumption.  Then we may inductively
define
$$\psi_{\cB,\cA}=\psi_{\cB(\lambda_N),\cB(\lambda_{N-1})}
\circ \ldots \circ \psi_{\cB(\lambda_1),\cB(\lambda_0)}.$$

Let $(C,s_1,\ldots,s_n)\stackrel{\pi}{\ra} B$
be a family of stable curves of type $(g,\cA)$.  
Under our simplifying assumption, we define 
$\psi_{\cB,\cA}(C,s_1,\ldots,s_n)$ as follows.  
Consider the ${\bQ}$-divisors $K_{\pi}(\cA):=K_{\pi}+a_1s_1+\ldots+a_ns_n$,
$K_{\pi}(\cB):=K_{\pi}+b_1s_1+\ldots+b_ns_n$, and 
$K_{\pi}(\cB(\lambda)):=
\lambda K_{\pi}(\cA)+(1-\lambda)K_{\pi}(\cB)$ for 
$\lambda\in \bQ, 0\le \lambda \le 1$.
We claim this is ample for each $\lambda \ne 0$.  This
follows from Remark \ref{rem:typesnew}.
If $L$ is nef but not ample, then there either exist sections
with weight zero (on type II components) or sets of
sections with weights summing to one (on type I components),
both of which are excluded by our assumptions.

We now apply Corollary \ref{cor:section} and Theorem \ref{thm:logcan}  
to obtain a new family of pointed nodal curves
$\pi':C'\ra B$ with smooth sections $s'_{i_1},\ldots,s'_{i_r}$
corresponding to the nonzero weights of ${\hat \cB}$.  
We define $\psi_{\cB,\cA}(C,s_1,\ldots,s_n)$ to be
the family $\pi':C'\ra B$, with image sections $s'_{i_1},\ldots,s'_{i_r}$
and weights $b_{i_1},\ldots,b_{i_r}$,  
a family of weighted pointed curves with
fibers in $\ocM_{g,\cB}$.  The vanishing statement
in Proposition \ref{prop:series} guarantees 
our construction commutes with
base extension.  Thus we obtain a natural transformation of
moduli functors
$$\psi_{\cB,\cA}:\ocM_{g,\cA} \ra \ocM_{g,\cB}.$$

Assume that $\cB={\hat \cB}$, so that $\psi_{\cB,\cA}$
is interpretted as a reduction map $\rho_{\cB,\cA}$.  Then
$\rho_{\cB,\cA}$ is an isomorphism over the locus of configurations
of points, it is a birational morphism.
$\square$

Reduction satisfies the following compatibility condition:
\begin{prop}
Fix $g$ and let $\cA,\cB,$ and $\cC$ be collections of weight data
so that the reductions $\rho_{\cB,\cA}, \rho_{\cC,\cB},$
and $\rho_{\cC,\cA}$ are well defined.  Then
$$\rho_{\cC,\cA}=\rho_{\cC,\cB}\circ \rho_{\cB,\cA}.$$
\end{prop}
{\em proof:}  Since the maps are birational morphisms of nonsingular
varieties, it suffices
to check that the maps coincide set-theoretically.  $\square$

\subsection{Exceptional loci}
\label{subsect:exceptional}
The exceptional locus of the reduction morphism is easily computed using
Corollary \ref{cor:section}: 
\begin{prop}\label{prop:reddesc}
The reduction morphism $\rho_{\cB,\cA}$ contracts the boundary divisors
$$D_{I,J}:=\oM_{0,\cA_I'} \times 
\oM_{g,\cA_J'}\quad \cA_I'=(a_{i_1},\ldots,a_{i_r},1) \quad
\cA_J'=(a_{j_1},\ldots,a_{j_{n-r}},1)$$
corresponding to partitions
$$\{1,\ldots,n \}=I\cup J \quad I=\{i_1,\ldots,i_r\} \quad J=\{j_1,\ldots, j_{n-r}\}$$
with $b_I:=b_{i_1}+\ldots +b_{i_r} \le 1$ and $2<r\le n.$
We have a factorization of
$\rho_{\cB,\cA}|D_{I,J}$:
$$\oM_{0,\cA_I'} \times \oM_{g,\cA_J'}\stackrel{\pi}{\longrightarrow}\oM_{g,\cA_J'} 
\stackrel{\rho}{\longrightarrow} \oM_{g,\cB_J'}\quad
\cB_J'=(b_{j_1},\ldots,b_{j_{n-r}},b_I)$$
where $\rho=\rho_{\cB'_J,\cA'_J}$ and $\pi$ is the projection.
\end{prop}

\begin{rem}\label{rem:blowdown}
Consider a set of weights $(a_{i_1},\ldots,a_{i_r})$ so that
$$a_{i_1}+\ldots +a_{i_r}> 1$$
but any proper subset has sum at most one.
Then $\oM_{0,\cA'_I}$
is isomorphic to ${\bP}^{r-2}$ (cf. \S \ref{subsect:Kapranov2}) and
$\rho_{\cB,\cA}$ is the blow-up of $\oM_{g,\cB}$ along the image
of $D_{I,J}$.  \end{rem}

\begin{cor}\label{cor:reddesc2}
Retain the notation and assumptions of Proposition \ref{prop:reddesc}.  
Assume in addition that 
for each $I\subset \{1,\ldots,n\}$ such that
$$a_{i_1}+\ldots +a_{i_r}> 1 \text{ and } b_{i_1}+\ldots +b_{i_r}\le 1$$
we have $r=2$.  Then $\rho_{\cA,\cB}$ is an isomorphism.
\end{cor}

\section{Chambers and walls}
Let $\cD_{g,n}$ denote the domain of admissible weight data
$$\cD_{g,n}:=\{
(a_1,\ldots,a_n)\in \bR^n :  0<a_j\le 1
\text{ and }
a_1+a_2+\ldots +a_n>2-2g \}.$$
A {\em chamber decomposition} of $\cD_{g,n}$ consists of
a finite set $\cW$ of hyperplanes $w_S\subset \cD_{g,n}$, the {\em walls}
of the chamber decomposition;  
the connected components of the complement to the
union of the walls $\cup_{S\in \cW}w_S$
are called the {\em open chambers.}

There are two natural chamber decompositions for $\cD_{g,n}$.
The {\em coarse chamber decomposition} is
$$\cW_c=\{ \sum_{j\in S}a_j=1: S\subset \{1,\ldots,n \}, 
2<|S|<n-2 \}$$
and the {\em fine chamber decomposition} is
$$\cW_f=\{ \sum_{j\in S}a_j=1: S\subset \{1,\ldots,n \},
2\le |S|\le n-2 \}.$$
Given a nonempty wall $w_S$, the set $\cD_{g,n}\setminus w$
has two connected components defined by the inequalities
$\sum_{j\in S}a_j>1$ and $\sum_{j\in S}a_j<1$.  

\begin{prop}\label{prop:chamberconst}
The coarse chamber decomposition is the coarsest decomposition
of $\cD_{g,n}$ such that $\ocM_{g,\cA}$ is constant on each chamber.  
The fine chamber decomposition is the coarsest decomposition
of $\cD_{g,n}$ such that $\cC_{g,\cA}$ is constant on each chamber.  
\end{prop}
{\em proof:}  It is clear that $\ocM_{g,\cA}$ (resp. $\cC_{g,\cA}$)
changes as we pass from one coarse (resp. fine) chamber to another.
It suffices then to show that $\ocM_{g,\cA}$ (resp. $\cC_{g,\cA}$)
is constant on each chamber.

Let $\cA$ and $\cB$ be contained in the
interior of a fine chamber and let  
$\pi:(C,s_1,\ldots,s_n) \ra B$
be a family in $\ocM_{g,\cA}(B)$.  Repeating the argument for
Theorem \ref{thm:reduction}, we find that
$K_{\pi}+b_1s_1+\ldots+b_ns_n$
is ample.  An application of Theorem \ref{thm:logcan} implies
that $(C,s_1,\ldots,s_n)$ is also $\cB$-stable.  Thus we get
an identification of $\ocM_{g,\cA}(B)$ and $\ocM_{g,\cB}(B)$
so $\cC_{g,\cA}\simeq \cC_{g,\cB}$.  

Consider the fine chambers contained
in a given coarse chamber $Ch$.
We shall say that two such fine chambers are {\em adjacent} 
if there exists a wall $w$ which is a codimension-one face of each.
Any two fine chambers in $Ch$ are related by a finite
sequence of adjacencies, so it suffices to show that adjacent
fine chambers yield the same moduli space.  
Fix two fine chambers in $Ch$ adjacent along $w$, which
we may assume is defined by $a_1+a_2=1$.  Let
$(w_1,\ldots,w_n)\in w$ be an element contained in the closure
of the chambers but not in any walls besides $w$.  For small
$\epsilon>0$, the weight data
$$\cA=(w_1+\epsilon,w_2+\epsilon,w_3,\ldots,w_n) \quad 
\cB=(w_1-\epsilon,w_2-\epsilon,w_3,\ldots,w_n)$$
lie in our two fine chambers.  Corollary \ref{cor:reddesc2}
implies that $\rho_{\cA,\cB}$ is an automorphism. $\square$

\begin{prob}
Find formulas for the number of nonempty walls
and chambers for $\cD_{g,n}$. 
\end{prob}

\begin{prop} \label{prop:perturb}
Let $\cA$ be a collection of weight data.  There exists
a collection of weight data $\cB$, contained in a fine open
chamber, such that $\cC_{g,\cA}=\cC_{g,\cB}$.
\end{prop}  
The fine open chamber produced in Proposition \ref{prop:perturb}
is called the fine chamber associated to $\cA$.

\noindent {\em proof:}  Consider the collection $T_1$
(resp. $T_{< 1}$, resp. $T_{> 1}$) of all subsets
$S\subset \{1,\ldots,n \}$ with $2\le |S| \le n-2$
such that $\sum_{j\in S}a_j=1$
(resp. $\sum_{j\in S}a_j<1$, resp. $\sum_{j\in S}a_j>1$.)
There exists a positive constant $\epsilon$ such that
for each $S\in T_{< 1}$ we have $\sum_{j\in S}a_j<1-\epsilon$,
for each $S\in T_{> 1}$ we have $\sum_{j\in S}a_j>1+\epsilon$,
$a_1+\ldots+a_n>2+\epsilon$, and $a_n>\epsilon$.  Setting 
$$\cB=(a_1-\epsilon/n,a_2-\epsilon/n,
\ldots, a_n-\epsilon/n)$$
and using results of \S \ref{subsect:exceptional},
we obtain the desired result.  $\square$

\begin{prop}
Let $\cA$ be weight data contained in a fine open
chamber.  Then $\cC_{g,\cA}$ is isomorphic to 
$\cM_{g,\cA\cup \{ \epsilon \}}$ for some sufficiently small
$\epsilon>0$.
\end{prop}
{\em proof:}  We retain the notation of the proof of
Proposition \ref{prop:perturb}.  Under our assumptions
the set $T_1$ is empty, so for each $T\subset \{1,\ldots,n\}$
with $2\le |S|\le n-2$ we have
$$|\sum_{j\in S}a_j -1 |>\epsilon. \quad \square$$

\section{Examples}
\label{sect:examples}

\subsection{Kapranov's approach to $\oM_{0,n}$}
\label{subsect:Kapranov}
The key classical observation
is that through each set of $n$ points of ${\bP}^{n-3}$
in linear general position, there passes a unique rational normal
curve of degree $n-3$.  It is therefore natural to realize
elements of $M_{0,n}$ as rational normal curves in projective
space.  This motivates Kapranov's blow-up construction of $\oM_{0,n}$
(\cite{Kap1} \S 4.3, \cite{Kap2}).
        
Start with $W_{1,1}[n]:={\bP}^{n-3}$ and choose $n-1$ points
$q_1,\ldots,q_{n-1}$ in linear general position.  
These are unique up to a projectivity.  We blow up as follows:
\begin{enumerate}
\item[ 1:]blow up the points $q_1,\ldots,q_{n-2}$, then the lines
passing through pairs of these points, followed by the planes
passing through triples of triples of these points, etc.;  

\item[ 2:]blow up the point $q_{n-1}$, then the lines spanned by
the pairs of points including $q_{n-1}$ but not $q_{n-2}$, then
the planes spanned by triples including $q_{n-1}$ but not
$q_{n-2}$, etc.; \ldots 
\item[ r:]blow up the linear spaces spanned by subsets 
$$\{q_{n-1},q_{n-2}\ldots,q_{n-r+1}\}\subset S \subset 
\{q_1,\ldots,q_{n-r-1},q_{n-r+1},\ldots,q_{n-1}\}$$
so that the order of the blow-ups
is compatible with the partial order on the subsets by inclusion;\ldots
\end{enumerate} 
Let $W_{r,1}[n],\ldots,W_{r,n-r-2}[n](:=W_r[n])$ denote the 
varieties produced at the $r$th step.  Precisely,
$W_{r,s}[n]$ is obtained once we finish blowing up subspaces
spanned by subsets $S$ with $|S|\le s+r-2$.    
Kapranov proves $\oM_{0,n}\simeq W_{n-3}[n]$.  

This may be interpretted with the reduction formalism 
of \S \ref{sect:reduction}.  
The exceptional divisors of the blow-downs $W_{r,s}[n]\ra W_{r,s-1}[n]$ are 
proper transforms of the boundary divisors 
$D_{I,J}$ corresponding to partitions
$$\{1,\ldots,n\}=I\cap J \text{ where } J=\{n\} \cup S.$$
Using the weight data
\begin{eqnarray*} 
\cA_{r,s}[n]&:=&
(\underbrace{1/(n-r-1),\ldots,1/(n-r-1)}_{n-r-1 \text{ times }},s/(n-r-1),
\underbrace{1,\ldots,1}_{r \text{ times }}), \\
& &r=1,\ldots,n-3,\quad s=1,\ldots,n-r-2,
\end{eqnarray*}
we realize $W_{r,s}[n]$ as $\oM_{0,\cA_{r,s}[n]}$.  
The blow-downs defined above are the reduction morphisms
\begin{eqnarray*}
\rho_{\cA_{r,s-1}[n],\cA_{r,s}[n]}&:&\oM_{0,\cA_{r,s}[n]} \longrightarrow \oM_{0,\cA_{r,s-1}[n]}
\quad s=2,\ldots,n-r-2\\
\rho_{\cA_{r,n-r-2}[n],\cA_{r+1,1}[n]}&:&\oM_{0,\cA_{r+1,1}[n]} \longrightarrow \oM_{0,\cA_{r,n-r-2}[n]}.
\end{eqnarray*}
In particular, Kapranov's factorization of $\oM_{0,n}\ra {\bP}^{n-3}$ as a sequence
of blow-downs is naturally a composition of reduction morphisms.

\subsection{An alternate approach to Kapranov's moduli space}
\label{subsect:Kapranov2}
There are many factorizations of the map $\oM_{0,n} \ra {\bP}^{n-3}$
as a composition of reductions.  We give another example here.
First, we compute the fine chamber containing the weight data
$\cA_{1,1}[n]$ introduced in \S \ref{subsect:Kapranov}:   
\begin{enumerate}
\item[]{$a_1+\ldots+{\widehat a_i}+\ldots+a_{n-1}\le 1$ 
(and thus $a_i+a_n>1$) for $i=1,\ldots,n-1;$}
\item[]{$a_1+\ldots+a_{n-1}>1$.}
\end{enumerate}
The corresponding moduli space $\oM_{0,\cA}$ is denoted $X_0[n]$
and is isomorphic to ${\bP}^{n-3}$.   
Similarly, we define $X_k[n]:=\oM_{0,\cA}$ provided $\cA$
satisfies
\begin{enumerate}
\item[]{$a_i+a_n>1$ for $i=1,\ldots,n-1;$}
\item[]{$a_{l_1}+\ldots+a_{l_r} \le 1$ (resp. $>1$) 
for each $ \{l_1,\ldots,l_r\}
\subset \{1,\ldots,n-1\}$ with $r\le n-k-2$ (resp. $r>n-k-2$).}
\end{enumerate}  
When $a_n=1$ and $a_1=\ldots=a_{n-1}=a(n,k)$
we have
$$1/(n-1-k)<a(n,k)\le 1/(n-2-k) $$
and 
$$\cA(n,k)=(\underbrace{a(n,k),\ldots,a(n,k)}_{n-1 \text{ times }},1).$$
Thus there exist reduction maps
$$\begin{array}{ccccc}
\rho_{\cA(n,k-1),\cA(n,k)}: & \oM_{0,\cA(n,k)} &\longrightarrow& 
\oM_{0,\cA(n,k-1)}& \\
 & || & & || & k=1,\ldots,n-4\\
 & X_k[n] &\longrightarrow & X_{k-1}[n]& 
\end{array}$$

We interpret the exceptional divisors of the induced maps.
For each partition
$$\{1,\ldots,n\}=I\cup J \quad I=\{i_1=n,i_2,\ldots,i_r \}, 
J=\{j_1,\ldots,j_{n-r}\} \ 2\le r \le n-2$$
consider the corresponding boundary divisor in $\oM_{0,n}$
$$D_{I,J}\simeq \oM_{0,r+1} \times \oM_{0,n-r+1}.$$
The divisors with $|I|=r<n-2$ are the exceptional divisors of
$X_{r-1}[n] \ra X_{r-2}[n]$.  The blow-down
$X_1[n] \ra X_0[n]$ maps $D_{I,J}$ with $I=\{n,i\}$ 
to the distinguished point
$q_i\in {\bP}^{n-3}$ mentioned in \S \ref{subsect:Kapranov}.  
Applying Proposition \ref{prop:reddesc} and 
$$(\underbrace{a(n,1),\ldots,a(n,1)}_{n-2\text{ times}},1)
\simeq (\underbrace{a(n-1,0),\ldots,a(n-1,0)}_{n-2 \text{ times }},1)
=\cA(n-1,0),$$
we obtain that the fiber over $q_i$ is 
$\oM_{0,\cA(n-1,0)}\simeq \bP^{n-4}$.
More generally, $\oM_{0,n}=X_{n-4}[n] \ra X_0[n]$
maps $D_{I,J}$ to the linear subspace spanned by 
$q_{i_2},\ldots,q_{i_r}$;  the fibers are isomorphic to ${\bP}^{n-r-2}$.   
The divisors $D_{I,J}$
with $|I|=n-2$ are proper transforms of
the hyperplanes of $X_0[n]\simeq {\bP}^{n-3}$ spanned by 
$q_{i_2},\ldots,q_{i_{n-2}}$.  

Thus our reduction maps give the following 
realization of $\oM_{0,n}$ as a blow-up of ${\bP}^{n-3}$:
\begin{enumerate}
\item[ 1:]{blow up the points $q_1,\ldots,q_{n-1}$
to obtain $X_1[n]$;}
\item[ 2:]{blow up proper transforms of
lines spanned by pairs of the points $q_1,\ldots,q_{n-1}$
to obtain $X_2[n]$;}
\item[ 3:]{blow up proper transforms of $2$-planes 
spanned by triples of the points to obtain $X_3[n]$;   \ldots}
\item[ n-4:]{blow up proper transforms of $(n-5)$-planes
spanned by $(n-4)$-tuples of the points to obtain $X_{n-4}[n]$.}
\end{enumerate}

\subsection{Keel's approach to $\oM_{0,n}$}
\label{subsect:Keel}
Let $U\subset ({\bP}^1)^n$ denote the set of all configurations
of $n$ {\em distinct} points in ${\bP}^1$, and $Q$ the resulting
quotient under the diagonal action of $\PGL_2$.  For each configuration
$(p_1,\ldots,p_n)$, there exists a unique projectivity $\phi$
with
$$\phi:(p_1,p_2,p_3) \ra (0,1,\infty).$$
The image of the configuration in $Q$ is determined completely by
the points $(\phi(p_4),\ldots,\phi(p_n))$ and we obtain an imbedding
$$Q\hookrightarrow ({\bP}^1)^{n-3}.$$ 

We may interpret $({\bP}^1)^{n-3}$ as $\oM_{0,\cA}$ where $\cA=(a_1,\ldots,a_n)$
satisfies the following inequalities:
\begin{enumerate}
\item[]{$a_{i_1}+a_{i_2}>1$ where $\{i_1,i_2\}\subset \{1,2,3\}$;}
\item[]{$a_i+a_{j_1}+\ldots+a_{j_r}\le 1$ for $i=1,2,3$ and
$\{j_1,\ldots,j_r\} \subset \{4,5,\ldots,n \}$ with $r>0$.}
\end{enumerate}
These conditions guarantee
that none of the first three sections coincide, but any
of the subsequent sections may coincide with any of the first three
or with one another.  When $a_1=a_2=a_3=a$ and
$a_4=a_5=\ldots=a_n=\epsilon$ we have 
$1/2 < a \le 1$ and $0<(n-3)\epsilon \le 1-a$.  Taking
$$\cA(n)=(a,a,a,\underbrace{\epsilon,\ldots,\epsilon}_{n-3 \text{ times }}),$$
we obtain that the compactification of $Q$ by $({\bP}^1)^{n-3}$ 
is isomorphic to $\oM_{0,\cA(n)}$.

We factor $\rho:\oM_{0,n} \ra (\bP^1)^{n-3}$ as a product
of reduction morphisms.  Let $\Delta_d$ denote the union of the
dimension $d$ diagonals, i.e., the locus where at least $n-2-d$
of the points coincide.  We will use this notation for both the locus
in $({\bP}^1)^{n-3}$ and its proper transforms.  Let
$$F_0=\pi_1^{-1}(0)\cup \pi_2^{-1}(0)\cup \ldots \cup \pi_{n-3}^{-1}(0)$$
be the locus of points mapping to $0$ under one of the projections $\pi_j$;
we define $F_1$ and $F_{\infty}$ analogously.  Again, we use the
same notation for proper transforms.  

Write $Y_0[n]=({\bP}^1)^{n-3}$ and define the first sequence
of blow-ups as follows:
\begin{enumerate}
\item[ 1:]{$Y_1[n]$ is the blow-up along the 
intersection $\Delta_1\cap (F_0 \cup F_1 \cup F_{\infty})$;}
\item[ 2:]{$Y_2[n]$ is the blow-up along the 
intersection $\Delta_2\cap (F_0 \cup F_1 \cup F_{\infty})$; \ldots}
\item[ n-4:]{$Y_{n-4}[n]$ is the blow-up along the 
intersection $\Delta_{n-4}\cap (F_0 \cup F_1 \cup F_{\infty})$.}
\end{enumerate}
The variety $Y_k[n]$ is realized by $\oM_{0,\cA}$ where
\begin{enumerate} 
\item[]{$a_{i_1}+a_{i_2}>1$ where $\{i_1,i_2\}\subset \{1,2,3\}$;}
\item[]{$a_i+a_{j_1}+\ldots+a_{j_r}\le 1$ (resp. $>1$) for $i=1,2,3$ and
$\{j_1,\ldots,j_r\} \subset \{4,5,\ldots,n \}$ with $0<r\le n-3-k$
(resp. $r>n-3-k$).}
\end{enumerate}
The second sequence of blow-ups is
\begin{enumerate}
\item[ n-3:]{$Y_{n-3}[n]$ is the blow-up along $\Delta_1$;}
\item[ n-2:]{$Y_{n-2}[n]$ is the blow-up along $\Delta_2$; \ldots}
\item[ 2n-9:]{$Y_{2n-9}[n]$ is the blow-up along $\Delta_{n-5}$.}
\end{enumerate}
The variety $Y_{k+n-4}[n]$ is realized by $\oM_{0,\cA}$ where
\begin{enumerate} 
\item[]{$a_{i_1}+a_{i_2}>1$ where $\{i_1,i_2\}\subset \{1,2,3\}$;}
\item[]{$a_{j_1}+\ldots+a_{j_r}\le 1$ (resp. $>1$) for
$\{j_1,\ldots,j_r\} \subset \{4,5,\ldots,n \}$ with $0<r\le n-3-k$
(resp. $r>n-3-k$).}
\end{enumerate}
\begin{rem}
Keel \cite{Ke} has factored $\rho$
as a sequence of blow-ups along smooth codimension-two centers
in the course of computing the Chow groups of $\oM_{0,n}$.  
However, the intermediate steps of his factorization do not
admit interpretations as $\oM_{0,\cA}$.  
For example, consider
$$\oM_{0,6}\ra \oM_{0,4}\times \oM_{0,5} \ra 
\oM_{0,4}\times \oM_{0,4}\times \oM_{0,4}\simeq ({\bP}^1)^3,$$
where the maps $\oM_{0,N} \ra \oM_{0,4}\times \oM_{0,N-1}$
are products of the forgetting morphisms arising from the subsets
$$\{1,2,3,N\},\{1,2,3,\ldots,N-1\}\subset \{1,2,\ldots,N\}.$$
The intermediate space $\oM_{0,4}\times \oM_{0,5}$ is not
of the form $\oM_{0,\cA}$ for any $\cA$.  
\end{rem}

\subsection{Losev-Manin moduli spaces}
\label{subsect:LM}
We refer to the paper \cite{LM}, where the following
generalization of stable pointed curves is defined.
This space was also studied by Kapranov (\cite{Kap1}, \S 4.3) 
as the closure of a generic orbit of $({\bC^*})^{n-3}$
in the space of complete flags in ${\bC}^{n-2}$.  

Let $S$ and $T$ be two finite disjoint sets
with $|S|=r$ and $|T|=n-r$, $B$ a scheme,
and $g\ge 0$.  An $(S,T)$-pointed stable curve of genus $g$ over
$B$ consists of the data:
\begin{enumerate}
\item{a flat family $\pi:C \ra B$ of nodal geometrically connected 
curves of arithmetic genus $g$;}
\item{sections $s_1,\ldots,s_r,t_1,\ldots,t_{n-r}$ 
of $\pi$ contained in the
smooth locus of $\pi$;}
\end{enumerate}
satisfying the following stability conditions:
\begin{enumerate}
\item{$K_{\pi}+\supp(s_1+\ldots +s_r+t_1+\ldots +t_{n-r})$ is $\pi$-relatively
ample;}
\item{each of the sections $s_1,\ldots,s_r$ are disjoint from all the other
sections, but $t_1,\ldots,t_{n-r}$ may coincide.}
\end{enumerate}
Now assume that $r=2$ and $g=0$, and consider
only pointed curves satisfying the following properties:
\begin{enumerate}
\item{the dual graph is linear;}
\item{the sections $s_1$ and $s_2$ are
contained in components corresponding to
the endpoints of the graph.}
\end{enumerate}
Then there is a smooth, separated, irreducible,
proper moduli space $\oL_{n-2}$ representing such
$(S,T)$-pointed stable curves (Theorem 2.2 of \cite{LM}).  

This space is isomorphic to $\oM_{0,\cA}$ where $\cA$ satisfies
the following conditions
\begin{enumerate}
\item[]{$a_1+a_i>1$ and $a_2+a_i>1$ for each $i=1,\ldots,n$;}
\item[]{$a_{j_1}+\ldots+a_{j_r}\le 1$ for each
$\{j_1,\ldots,j_r\} \subset \{3,\ldots,n \}$ with $r>0$.}
\end{enumerate}
We emphasize that these conditions force the dual graphs of
the associated curves to have the properties postulated
by Losev and Manin.  Specializing the weights, we
obtain that $\oL_{n-2}\simeq \ocM_{0,\cA}$ where    
$$\cA=(1,1,\underbrace{1/(n-2),\ldots,1/(n-2)}_{n-2 \text{ times }}).$$

Using the reduction maps we obtain explicit blow-up realizations
of the Losev-Manin moduli spaces.  
Setting $\cB=(1,1/(n-2),\ldots,1/(n-2))$
we obtain a morphism
$$\rho_{\cB,\cA}:\oM_{0,\cA} \ra \oM_{0,\cB}$$    
where $\oM_{0,\cB}\simeq\bP^{n-3}$ is Kapranov's compactification.
The Losev-Manin moduli space is
the first step of the factorization
described in \S \ref{subsect:Kapranov}.

\section{Interpretations as log minimal models of moduli spaces}
\label{sect:LMMP}

Recently, Keel and his collaborators 
\cite{KeMc1},\cite{GKM},\cite{G},\cite{Ru}
have undertaken a study of the birational geometry
of the moduli space of curves, with an emphasis on the geometry
of the cones of effective curves and divisors.  In many
cases, they find natural modular interpretations for 
contractions and modifications arising from the minimal model
program.  Therefore, one might expect that 
natural birational modifications of a moduli space $\oM$ should 
admit interpretations 
as log minimal models with respect to a boundary
supported on natural divisors of $\oM$.  The most accessible divisors
are those parametrizing degenerate curves, so we focus on these here.

Fix $\oM=\oM_{0,n}$, the moduli space of genus zero curves
with $n$ marked points.  For each unordered partition
$$\{1,2,\ldots,n\}=I\cup J \text{ where }|I|,|J|>1$$
let $D_{I,J}\subset \oM$ denote the divisor corresponding         
to the closure of the locus of pointed curves with two irreducible
components, with the sections indexed by $I$ on one component
and by $J$ on the other.  Let $\delta$
be the sum of these divisors, with each $D_{I,J}$ appearing
with multiplicity one.  We can also describe degenerate curves 
on $\oM_{0,\cA}$ (see \S \ref{subsect:exceptional}).  
There are two types
to consider.  First, consider a partition
as above satisfying
$$a_{i_1}+\ldots+a_{i_r}>1, a_{j_1}+\ldots+a_{j_{n-r}}>1\quad
I=\{i_1,\ldots,i_r \},  J=\{j_1,\ldots,j_{n-r}\}.$$
Let $D_{I,J}(\cA)$ denote the image of
$D_{I,J}$ in $\oM_{0,\cA}$ under the reduction map;  
this is a divisor and parametrizes nodal curves as above.
The union of such divisors is denoted $\nu$.    
Second, any partition
with $I=\{i_1,i_2 \}$ and $a_{i_1}+a_{i_2}\le 1$ also corresponds
to a divisor $D_{I,J}(\cA)$ in $\oM_{0,\cA}$, parametrizing
curves where the sections $s_{i_1}$ and $s_{i_2}$ coincide.  
These curves need not be nodal.  The union of such divisors
is denoted $\gamma$ and the union of $\gamma$ and $\nu$
is denoted $\delta$.  The remaining partitions yield
subvarieties of $\oM_{0,\cA}$ with codimension $>1$.

\begin{prob}\label{prob:LMMP}
Consider the moduli space $\oM_{0,\cA}$ of weighted pointed
stable curves of genus zero.  Do there exist
rational numbers $d_{I,J}$ so that
$$K_{\oM_{0,\cA}}+\sum_{I,J}d_{I,J}D_{I,J}(\cA)$$
is ample and log canonical?
The sum is taken over partitions
$$\{1,2,\ldots,n\}=I\cup J \text{ where }|I|,|J|\ge 2,$$
where either 
$$a_{i_1}+\ldots+a_{i_r}>1 \text{ and } a_{j_1}+\ldots+a_{j_{n-r}}>1$$
or 
$$r=2 \text{ and } a_{i_1}+a_{i_2}\le 1.$$ 
\end{prob}
The assertion that the singularities are log canonical
implies that the coefficients are nonnegative and $\le 1$.  

We shall verify the assertion of Problem \ref{prob:LMMP}
in examples by computing the discrepancies  
of the associated reduction morphisms.  We also refer
the reader to Remark \ref{rem:10sing} for another instance
where it is verified. 

\subsection{Mumford-Knudsen moduli spaces}
It is well known that 
$K_{\oM_{0,n}}+\delta$ 
is ample and log canonical.  We briefly sketch the 
proof.  For ampleness, we use the identity
$$\kappa_{\text{log}}:=
\pi_*[c_1(\omega_{\pi}(s_1+\ldots+s_n))^2]=
K_{\oM_{0,n}}+\delta$$
where 
$$\pi:(\cC_{0,n},s_1,\ldots,s_n)\ra \oM_{0,n}$$ 
is the universal curve.
Fix pointed elliptic curves $(E_i,p_i), i=1,\ldots,n$, which 
we attach to an $n$-pointed rational curve to obtain a stable curve
of genus $n$.  This yields an imbedding
$$j:\oM_{0,n} \ra \oM_n.$$
The divisor 
$$\kappa=u_*[c_1(\omega_u)^2],$$
where $u:\cC_n \ra \oM_n$ is the universal curve, is ample 
(\cite{HM} 3.110 and 6.40)
and pulls back to $\kappa_{\text{log}}$.
As for the singularity condition, it suffices to observe that
through each point of $\oM_{0,n}$ there pass at most $n-3$ boundary
divisors which intersect in normal crossings.

\subsection{Kapranov's examples}
We retain the notation of \S \ref{subsect:Kapranov2}.  
The boundary divisors on 
$$X_0[n]=\oM_{0,\cA}\simeq {\bP}^{n-3}$$ 
are indexed by partitions
$$\{1,2,\ldots,n\}=\{i_1,i_2 \} \cup \{j_1,\ldots,j_{n-3},n\}$$
and correspond to the hyperplanes spanned by the points 
$q_{j_1},\ldots,q_{j_{n-3}}$.  Consider the log canonical
divisor
$$K_{X_0[n]}+\alpha D[n], \text{ where } D[n]:=\sum_{I=\{i_1,i_2\}}D_{I,J}(\cA),$$
which is ample if and only if
\begin{equation}
\alpha \binom{n-1}{2}>n-2, \text{ i.e., } \alpha>2/(n-1).
\label{ineq1}
\end{equation}

Let $E_k[n]$ denote the exceptional divisor of 
the blow-up $X_k[n] \ra X_{k-1}[n]$ and 
$$\rho_k:X_k[n] \ra X_0[n]$$
the induced birational morphism.
Each component of $E_k[n]$ is obtained by blowing up a nonsingular
subvariety of codimension $n-2-k$, so we obtain
$$K_{X_k[n]}=\rho_k^*K_{X_0[n]} + \sum_{r=1}^k (n-3-r)E_r[n].$$
Through each component of the center of $E_k[n]$, there
are $\binom{n-1-k}{n-3-k}$ nonsingular
irreducible components of $D[n]$.  It follows that
$$\rho_k^*D[n]=D[n]_k+\sum_{r=1}^k\binom{n-1-r}{n-3-r}E_r[n]$$
where $D[n]_k$ is the proper transform of $D[n]$.  
The discrepancy equation takes the form
$$K_{X_k[n]}+\alpha D[n]_k=\rho^*_k(K_{X_0[n]}+\alpha D[n])+
\sum_{r=1}^k(n-3-r-\alpha\binom{n-1-r}{n-3-r})E_r[n],$$
which is log canonical provided
$$-1\ge (n-3-r)-\alpha\binom{n-1-r}{n-3-r} \quad r=1,\ldots,k.$$
This yields the condition
$$\alpha\le 2/(n-2),$$
which is compatible with inequality (\ref{ineq1}).  

This computation yields a positive answer
to Problem \ref{prob:LMMP} for $X_0[n-4]$ but not necessarily
for $X_k[n-4]$ with $k>0$.  The exceptional divisors $E_r[n]$
may have large positive discrepancies.
 
\subsection{Keel's example}
We retain the notation of \S \ref{subsect:Keel} 
and use $E_k[n]$ for the exceptional divisor of $Y_k[n]\ra Y_{k-1}[n]$
and $\rho_k$ for the birational morphism $Y_k[n] \ra Y_0[n]$.
Let $F[n]=F_0+F_1+F_{\infty}$ and $D[n]$ the union of the diagonals
in $({\bP}^1)^{n-3}$.  Their proper transforms are denoted
$F[n]_k$ and $D[n]_k$.  The divisor 
$K_{(\bP^1)^{n-3}}+\alpha F[n]+\beta D[n]$
is ample if and only if
\begin{equation}
3\alpha+(n-4)\beta > 2.\label{ineq2}
\end{equation}  
We have the following discrepancy equations
\begin{eqnarray*}
K_{Y_{2n-9}[n]}&=&\rho^*_{2n-9}K_{Y_0[n]}
	+\sum_{r=1}^{n-4}(n-3-r)E_r[n] + 
	\sum_{r=1}^{n-5}(n-4-r)E_{n-4+r}[n]\\
 \rho^*_{2n-9}F[n]&=& F[n]_{2n-9}+
	\sum_{r=1}^{n-4}(n-2-r)E_r[n]\\
 \rho^*_{2n-9}D[n]&=& D[n]_{2n-9}+\sum_{r=1}^{n-4}\binom{n-2-r}{2}E_r[n]
 	+\sum_{r=1}^{n-5}\binom{n-2-r}{2}E_{n-4+r}[n],
\end{eqnarray*}
which yield inequalities
\begin{eqnarray*}
-1&\le& (n-3-r)-\alpha(n-2-r)-\beta\binom{n-2-r}{2},r=1,\ldots,n-4 \\
-1&\le& (n-4-r)-\beta\binom{n-2-r}{2}, r=1,\ldots,n-5.
\end{eqnarray*}
These in turn yield
$$\beta\le 2/(n-3) \quad \alpha+\beta((n-4)/2)\le 1.$$
To satisfy these conditions and inequality (\ref{ineq2}),
we may choose
$$\alpha=1/(n-3) \quad \beta=2/(n-3).$$

\section{Variations of GIT quotients of $(\bP^1)^n$}
\label{sect:GIT}
In this section, we show how geometric invariant
theory quotients of $(\bP^1)^n$ may be interpretted 
as `small-parameter limits' of the moduli
schemes $\oM_{0,\cA}$ as $\sum_{j=1}^n a_j\rightarrow 2$.

We review the description of the stable locus for 
the diagonal action of $\PGL_2$ on $(\bP^1)^n$
(see \cite{Th}, \S 6 and \cite{GIT}, Ch. 3).  The
group $\PGL_2$ admits no characters, so ample
fractional linearizations correspond to line bundles
$\cO(t_1,\ldots,t_n)$ on $({\bP}^1)^n$, where the $t_i$
are positive rational numbers.  
A point $(x_1,\ldots,x_n)\in (\bP^1)^n$ is {\em stable} 
(resp. {\em semistable}) if, for each $x\in {\bP}^1$,
$$\sum_{j=1}^nt_j\delta(x,x_j) < (\le)
\frac{1}{2}(\sum_{j=1}^n t_j),$$
where $\delta(x,x_j)=1$ when $x=x_j$ and $0$ otherwise.
\begin{rem}\label{rem:linopen}
The linearizations $\cT=(t_1,\ldots,t_n)$
for which a given point of $(\bP^1)^n$ is stable
(resp. semistable) form an open (resp. closed) subset.  
\end{rem}   

To strengthen the analogy with our moduli spaces, 
we renormalize so that $t_1+t_2+\ldots+ t_n=2$.  
Then the stability condition takes the following
form:  for any $\{i_1,\ldots,i_r \} \subset \{ 1,\ldots,n \}$,
$x_{i_1},\ldots,x_{i_r}$ may coincide only when
$$t_{i_1}+\ldots +t_{i_r} < 1.$$
In particular, the stable locus is nonempty only when each $t_j<1$.
In this case, the corresponding GIT quotient is denoted $\cQ(t_1,\ldots,t_n)$
or $\cQ(\cT)$.  We define the {\em boundary} of $\cD_{0,n}$ as
$$\partial\cD_{0,n}:=\{(t_1,\ldots,t_n): t_1+\ldots+t_n=2, 0<t_i<1 
\text{ for each } i=1,\ldots,n\}.$$ 
We shall say that $\cT$ is {\em typical} if all semistable points
are stable and {\em atypical} otherwise.
Of course, $\cT$ is typical exactly when $t_{i_1}+\ldots +t_{i_r}\ne 1$
for any $\{i_1,\ldots,i_r \}\subset \{1,\ldots,n \}$.  

\begin{thm}\label{thm:typical}
Let $\cT$ be a typical linearization in $\partial \cD_{0,n}$.
Then there exists an open neighborhood $U$ of $\cT$
so that $U\cap \cD_{0,n}$ is contained in an open
fine chamber of $\cD_{0,n}$.  For each set of weight data
$\cA\in U\cap \cD_{0,n}$, there is a natural isomorphism
$$\oM_{0,\cA} \stackrel{\simeq}{\longrightarrow} \cQ(\cT).$$
\end{thm}
{\em proof:}
We choose 
\begin{eqnarray*}
U&=&\{(u_1,\ldots,u_n)\in {\bQ}^n: 0<u_i<1 \text{ and } \\ 
& &u_{i_1}+\ldots+u_{i_r}\ne 1 \text{ for any }
\{i_1,\ldots,i_r\}\subset \{1,\ldots,n\}\}.
\end{eqnarray*}
It follows that $U\cap \cD_{0,n}$ is contained in a 
fine chamber.  In particular, Proposition 
\ref{prop:chamberconst} implies that for any weight data
$\cA_1,\cA_2\in U\cap \cD_{0,n}$, we have
$\oM_{0,\cA_1} \simeq \oM_{0,\cA_2}$.  

We construct the morphism $\cQ(\cT)\ra \oM_{0,\cA}$,
where $\cA\in U\cap \cD_{0,n}$.  
Points of $\cQ(\cT)$ classify $(x_1,\ldots,x_n)\in (\bP^1)^n$
up to projectivities, where $x_{i_1},\ldots,x_{i_r}$ do not
coincide unless $t_{i_1}+\ldots+t_{i_r}< 1$.  
In our situation, $t_{i_1}+\ldots + t_{i_r}< 1$
if and only if $a_{i_1}+\ldots +a_{i_r}< 1$.  Since
$K_{\bP^1}+\sum a_jx_j$ has positive degree,
we may conclude that $(\bP^1,x_1,\ldots,x_n)$
represents a point of $\oM_{0,\cA}$.  As
$\cQ(\cT)$ parametrizes a family of pointed curves in $\oM_{0,\cA}$,
we obtain a natural morphism $\cQ(\cT)\ra \oM_{0,\cA}$.  This is 
a bijective birational projective morphism from a normal variety to a regular variety,
and thus an isomorphism.  $\square$

\begin{figure}[h]
\centerline{\psfig{figure=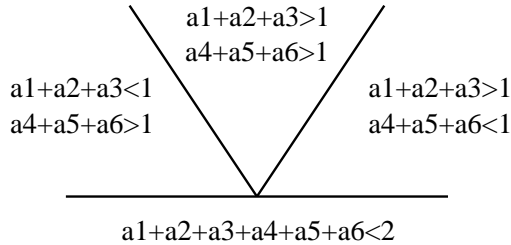}}
\caption{Chamber structure at the boundary}
\label{fig:weight1}
\end{figure}

For an atypical point $\cT$ of the boundary the description 
is more complicated.  From a modular standpoint, 
each neighborhood of $\cT$
intersects a number of fine chambers arising from 
different moduli problems.
See Figure \ref{fig:weight1} for 
a crude schematic diagram.
From an invariant-theoretic standpoint, for each
linearization $\cT'$ in a sufficiently small neighborhood
of $\cT$ we have a birational morphism
$$\cQ(\cT')\ra \cQ(\cT),$$
induced by the inclusion of the $\cT'$-semistable points
into the $\cT$-semistable points (cf. Remark
\ref{rem:linopen}).  A more general discussion of
this morphism may be found in \S 2 of \cite{Th}.    

\begin{thm}\label{thm:atypical}
Let $\cT\in \partial \cD_{0,n}$ be an 
atypical linearization.  Suppose that $\cT$ is
in the closure of
the coarse chamber associated with the weight data $\cA$.
Then there exists a natural
birational morphism
$$\rho: \oM_{0,\cA} \longrightarrow \cQ(\cT).$$
\end{thm}
{\em proof:}  
Consider the map $\tau: \cD'_{0,n} \ra \partial\cD_{0,n}$
given by the rule
$$\tau(b_1,\ldots,b_n)=(b_1B^{-1},\ldots,b_nB^{-1}),$$
where $B=(b_1+\ldots +b_n)/2$ and
$$\cD'_{0,n}=\{(b_1,\ldots,b_n):b_1+\ldots+b_n\ge 2,
0 < b_i < 1 \}.$$
For each $\cB\in \cD'_{0,n}$ such that $\tau(\cB)$ is typical, we
obtain a birational morphism
$$\rho_{\tau(\cB),\cB}:\oM_{0,\cB} \ra \cQ(\tau(\cB)).$$
This is defined as the composition of the reduction 
morphism $\rho_{\cB_1,\cB}:\oM_{0,\cB}\ra \oM_{0,\cB_1}$,
where $\cB_1=\epsilon\cB+(1-\epsilon)\tau(\cB)$ for small
$\epsilon>0$, and the isomorphism 
$\oM_{0,\cB_1} \ra \cQ(\tau(\cB))$ given by Theorem \ref{thm:typical}.        
 
The closure of each coarse chamber is the union of the closures
of the finite collection of fine chambers contained in it.  
By Proposition \ref{prop:chamberconst}, we may assume $\cA$ is 
in a fine open chamber $Ch$ with 
closure $Ch'\subset \cD'_{0,n}$ containing $\cT$.  
Clearly, $\tau(Ch')$ contains $\cT$ and $\tau(Ch)$ 
contains typical points arbitrarily close to $\cT$.  Choose
$\cA\in Ch$ so that $\tau(\cA)$ is typical and close to $\cT$, so 
there exists a generalized reduction morphism
$$\rho_{\tau(\cA),\cA}:\oM_{0,\cA} \ra \cQ(\tau(\cA))$$
and an induced birational morphism of GIT quotients
$$\cQ(\tau(\cA)) \ra \cQ(\cT).$$
Composing, we obtain the birational morphism claimed in the
theorem. $\square$

In light of Theorems \ref{thm:typical} and \ref{thm:atypical},
when $\cT \in \partial \cD_{0,n}$
we may reasonably interpret the GIT quotient $\cQ(\cT)$
as $\oM_{0,\cT}$.  This gives one possible definition for
moduli spaces with weights summing to two (cf. \S \ref{subsect:wtstwo}.)

\begin{rem} (Suggested by I. Dolgachev)
Theorem \ref{thm:typical} implies that
$\oM_{0,\cA}$ is realized as
a GIT quotient $Q(\cT)$ when the closure 
of the coarse chamber associated with 
$\cA$ contains a typical linearization.
For example, if
$$\cA=(2/3,2/3,2/3,
\underbrace{\epsilon,\ldots,\epsilon}_{n-3 \text{ times }}), \quad
\epsilon>0 \text{ small,}$$
then the moduli space $\oM_{0,\cA}\simeq (\bP^1)^{n-3}$,
studied in \S \ref{subsect:Keel},
arises as a GIT quotient (see \cite{KLW}).  
\end{rem}

\begin{rem}\label{rem:10sing}
We explicitly construct 
$$\oM_{0,\cT}:=\cQ(\cT)  \quad \cT=(1/3,1/3,1/3,1/3,1/3,1/3)$$
using a concrete description of the map
$$\rho: \oM_{0,(1,1/3,1/3,1/3,1/3,1/3)} 
\ra \oM_{0,\cT}$$
produced in Theorem \ref{thm:atypical}.  
Recall that the first space is the space $X_1[6]$ 
(see \S \ref{subsect:Kapranov2}), isomorphic to ${\bP}^3$
blown-up at five points $p_1,\ldots,p_5$
in linear general position.  The map $\rho$
is obtained by contracting the proper transforms $\ell_{ij}$ of the ten lines 
joining pairs of the points.  It follows that $\oM_{0,\cT}$
is singular at these ten points.  Concretely, $\rho$ 
is given by the linear series of quadrics
on ${\bP}^3$ passing through $p_1,\ldots,p_5$.  
Thus we obtain a realization
of $\oM_{0,\cT}$ as a cubic hypersurface in ${\bP}^4$
with ten ordinary double points.  Finally, we observe that 
$K_{\oM_{0,\cT}}+\alpha\delta$ is ample and log canonical provided
$2/5<\alpha\le 1/2$, thus yielding a positive answer
to Problem \ref{prob:LMMP} in this case.    
\end{rem}

\

\noindent 
The Institute of Mathematical Sciences,\\
Room 501, Mong Man Wai Building,\\
Chinese University of Hong Kong,\\
Shatin, Hong Kong

\

\noindent {\em Current address:} \\
Rice University, \\
Math Department--MS 136, \\
6100 S. Main St., \\
Houston, Texas 77005-1892, \\
hassett@math.rice.edu 

\end{document}